\documentclass[aos,authoryear]{imsart}

\usepackage{
amssymb
}

\usepackage{amsthm,amsmath,natbib}
\startlocaldefs
\newtheorem{theorem}{Theorem}  % [section]

\newtheorem{lemma}{Lemma}
\newtheorem{corollary}{Corollary}

\newtheorem{remark}{Remark}
\endlocaldefs

\begin{document}

\begin{frontmatter}

% "Title of the paper"
\title{Empirical risk minimization in inverse problems: 
Extended technical version}
\runtitle{Empirical risk minimization}

% indicate corresponding author with \corref{}
% \author{\fnms{John} \snm{Smith}\corref{}\ead[label=e1]{smith@foo.com}\thanksref{t1}}
% \thankstext{t1}{Thanks to somebody} 
% \address{line 1\\ line 2\\ printead{e1}}
% \affiliation{Some University}

\author{\fnms{Jussi} \snm{Klemel\"a}\corref{}\ead[label=e1]{klemela@oulu.fi}}
\address{
Jussi Klemel\"a\\
University of Oulu\\
Department of Mathematical Sciences\\
P. O. Box 3000\\
90014 University of Oulu\\ 
Finland\\
Fax: +358-8-5531730\\
\printead{e1}}
\and
\author{\fnms{Enno} \snm{Mammen}\ead[label=e2]{emammen@rumms.uni-mannheim.de}}
\address{
Enno Mammen\\
University of Mannheim,
Department of Economics\\
L7 3-5,\\
68131 Mannheim, Germany\\
Fax +49-621-1811931\\
\printead{e2}}
\affiliation{University of Oulu and University of  Mannheim}

\runauthor{}

\begin{abstract}
We study estimation of a multivariate function 
$f:{\bf R}^d \to {\bf R}$ 
when the observations are available from function $Af$, where
$A$ is a known linear operator. Both the Gaussian white noise model
and density estimation are studied. We define an $L_2$ empirical
risk functional, which is used to define an $\delta$-net minimizer
and a dense empirical risk minimizer. Upper bounds for the mean
integrated squared error of the estimators are given. The upper
bounds show how the difficulty of the estimation depends on the
operator through the norm of the adjoint of the inverse of the
operator, and on the underlying function class through the entropy
of the class. Corresponding lower bounds are also derived. As
examples we consider convolution operators and the Radon transform.
In these examples the estimators achieve the optimal rates of
convergence.
Furthermore, a new type of oracle inequality is given for
inverse problems in additive models.
\end{abstract}

\begin{keyword}[class=AMS]
\kwd[Primary ]{62G07}
%\kwd{}
%\kwd[; secondary ]{}
\end{keyword}

\begin{keyword}
\kwd{deconvolution}
\kwd{empirical risk minimization}
\kwd{multivariate density estimation}
\kwd{nonparametric function estimation}
\kwd{Radon transform}
\kwd{tomography}
\end{keyword}

\end{frontmatter}

\section{Introduction}

We consider estimation of a function $f:{\bf R}^d \to {\bf R}$, when
a linear transform $Af$ of the function is observed under stochastic
noise. We consider both the Gaussian white noise model and density
estimation with i.i.d.~observations. We study two estimators: a
$\delta$-net estimator which minimizes the $L_2$ empirical risk over
a minimal $\delta$-net of a function class, and a dense empirical
risk minimizer which minimizes the empirical risk over the whole 
function class without
restricting the minimization over a $\delta$-net. We call this
estimator ``dense minimizer'' because it is defined as a minimizer
over a possibly uncountable function class. The $\delta$-net
estimator is more universal: it may be applied also for unsmooth
functions and for severely ill-posed operators. On the other hand,
the dense empirical minimizer is expected to work only for
relatively smooth cases (the entropy integral has to converge). But
because the minimization in the calculation of this estimator is not
restricted to a $\delta$-net we have available a larger toolbox of
algorithms for finding (an approximation of) the minimizer of the empirical
risk.

Let $({\bf Y},{\cal Y},\nu)$ be a Borel space
and
let $A: L_2({\bf R}^d) \to L_2({\bf Y})$ be a linear operator,
where 
$L_2({\bf R}^d)$ is the space of square integrable functions 
$f:{\bf R}^d\to {\bf R}$
(with respect to the Lebesgue measure), and
$L_2({\bf Y})$ is the space of square integrable functions 
$g: {\bf Y} \to {\bf R}$ (with respect to measure $\nu$).
In the density estimation model
we have i.i.d.~observations
\begin{equation}
Y_1,\ldots ,Y_n \in {\bf Y},
\end{equation}
with common density function $Af :{\bf Y} \to {\bf R}$, where
$f:{\bf R}^d \to {\bf R}$ is a density function which we want to
estimate.
In the Gaussian white noise model the observation is a realization of the
process
\begin{equation}  \label{yn}
dY_n(y) = (Af)(y)\,dy + n^{-1/2} \, dW(y),
\qquad
y \in {\bf Y},
\end{equation}
where
$W(y)$ is the Brownian process on ${\bf Y}$,  %${\bf R}^d$.
that is, for $h_1,h_2 \in L_2({\bf Y})$, the random vector
$(\int_{\bf Y} h_1dW, \int_{\bf Y} h_2dW)$ is a 2-dimensional Gaussian random
vector with $0$ mean, marginal variances
$\|h_1\|_{2,\nu}^2,\|h_2\|_{2,\nu}^2$, and covariance $\int_{\bf Y}
h_1h_2\,d\nu$. 
(In our examples ${\bf Y}$ is either the Euclidean space or the product
of the real line with the unit sphere, so that the existence of the
Brownian process is guaranteed.)
We want to estimate the signal function $f:{\bf R}^d
\to {\bf R}$. 
The Gaussian white noise model is very useful in
presenting the basic mathematical ideas in a transparent way. For
the $\delta$-net estimator the treatment is almost identical for the
Gaussian white noise model and for the density estimation, but when
we consider the dense empirical risk minimization, then in the density 
estimation model we need to
use  bracketing numbers and empirical
entropies with bracketing, instead of the usual $L_2$ entropies. 
Our results for the
Gaussian white noise model can also serve as first step for getting
analogous results for inverse problems in regression or in other
statistical models.

The $L_2$ empirical risk is defined by
\begin{equation} \label{gamman}
\gamma_n(g)
=
\left\{ \begin{array}{ll}
-2 \int_{{\bf Y}} (Qg) \,dY_n + \left\| g \right\|_2^2 ,
& \mbox{ Gaussian white noise, }
\\
-2 n^{-1} \sum_{i=1}^n (Qg)(Y_i) + \left\| g \right\|_2^2 ,
& \mbox{ density estimation, }
\end{array} \right.
\end{equation}
where $Q$ is the adjoint of the inverse of $A$:
\begin{equation}  \label{q-def}
\int_{{\bf R}^d}  (A^{-1}h)g
=
\int_{{\bf Y}}  h(Qg)\, d\nu ,
\end{equation}
for $h\in L_2({\bf Y})$, $g \in L_2({\bf R}^d)$.
The operator $Q = (A^{-1})^*$ has the domain $L_2({\bf R}^d)$,
similarly as $A$.
Minimizing $\| \hat{f}-f\|_2^2$ with respect to estimators $\hat{f}$ is
equivalent to minimizing
$\| \hat{f}-f\|_2^2 - \| f\|_2^2$,
and we have,
in the Gaussian white noise model,
\begin{eqnarray}
\| \hat{f}-f\|_2^2 - \| f\|_2^2
&=&  \nonumber
-2 \int_{{\bf R}^d} f\hat{f} + \|\hat{f}\|_2^2
\\ &=&  \nonumber
-2 \int_{{\bf Y}} (Af)(Q\hat{f}) \, d\nu + \|\hat{f}\|_2^2
\\ &\approx&  \nonumber
-2 \int_{{\bf Y}} (Q\hat{f}) \, dY_n
%-\, \frac{2}{n} \sum_{i=1}^n (Q\hat{f})(Y_i)
+ \|\hat{f}\|_2^2
\\ &=&  \label{empi-moti}
\gamma_n\left( \hat{f} \right) .
\end{eqnarray}
The usual least squares estimator is defined as a minimizer of the
the criterion
\begin{eqnarray}
\| A\hat{f} - Af \|_{\bf Y}^2 - \|Af\|_{\bf Y}^2
& \approx &  \nonumber
-2 \int_{{\bf Y}} (Ag) \,dY_n + \left\|Ag\right\|_{\bf Y}^2
\\ & \stackrel{def}{=} &
\tilde{\gamma}_n(g)  .
\label{least-risk}
\end{eqnarray}
See for example \cite{OSull86}.
In density estimation the log-likelihood empirical risk
has been more common than the $L_2$ empirical risk,
and in the setting of inverse problems the log-likelihood
is defined as
$\bar{\gamma}_n (g) = - n^{-1} \sum_{i=1}^n \log (Ag)(Y_i)$,
analogously to (\ref{least-risk}).
These alternative definitions of the empirical risk do not seem to lead to such
elegant theory as the empirical risk in (\ref{gamman}).
The empirical risk in (\ref{gamman}) has been used
in deconvolution problems for projection estimators by
\cite{ComteTaupiRozen05}.

We give upper bounds for the mean integrated squared error (MISE)
of the estimators.
The upper bounds characterize how
the rates of convergence depend on the entropy of the underlying function class
${\cal F}$ and on smoothness properties of the operator $A$.
Previously such characterizations have been given
(up to our knowledge) in inverse problems
only for the case of estimating real valued linear functionals $L$.
In these cases the rates of convergence are determined
by the modulus of continuity of the functional
$\omega(\epsilon) = \sup \{ L(f) : f\in {\cal F} , \|Af\|_2 \leq \epsilon\}$,
see \cite{DonohLow92}.
For the case of estimating the whole function with a global loss function
the rates of convergence
depend on the largeness of the underlying function class in terms of the
entropy and capacity,
see
\cite{Cenco72},
\cite{LeCam73},
\cite{IbragHasmi80},
\cite{IbragHasmi81},
\cite{Birge83},
\cite{HasmiIbrag90},
\cite{BarroYang99},
\cite{Ibrag04}.
$\delta$-net estimators were considered e.g.\ by
\cite{LaanDudoiVaart04}.
 These papers consider direct statistical problems.
We show that for inverse statistical problems the rate of convergence
depends on the operator trough the operator norm
$\varrho(Q,{\cal F}_{\delta})$ of $Q$,
over a minimal $\delta$-net ${\cal F}_{\delta}$, see (\ref{opnormQ})
for the definition of $\varrho(Q,{\cal F}_{\delta})$.
More precisely, the convergence rate $\psi_n$ of the $\delta$-net estimator
is the solution to the equation
$$
n\psi_n^2 = \varrho^2(Q,{\cal F}_{\psi_n}) \log (\# {\cal F}_{\psi_n}) ,
$$
where
$\# {\cal F}_{\psi_n}$ is the cardinality of a minimal $\delta$-net.
For direct problems, when $A$ is the identity operator,
$\varrho(Q,{\cal F}_{\delta}) \asymp 1$.
As examples of operators $A$ we consider the convolution operator
and the Radon transform.
For these operators the estimators achieve the minimax rates of convergence
over Sobolev classes.

The general framework for empirical risk minimization
and the use of the empirical process machinery including entropy bounds for
deriving optimal bounds seems to be new.
Convolution and Radon transforms are discussed for illustrative purposes.
These examples show that our results lead to optimal rates of convergence.
As a new application we introduce  
the estimation of additive models in inverse problems.
A new type of oracle inequality is presented, 
which gives the optimal rates of convergence also
in ``anisotropic'' inverse problems.

\paragraph{Contents}
Section~\ref{deltanetsec} gives an upper bound for the MISE of the
$\delta$-net estimator.
Section~\ref{lower-bound-section} gives a lower bound for the
MISE of any estimator.
Section~\ref{deepsec} gives an upper bound for the MISE
of the dense empirical risk minimizer.
Section~\ref{opesec} finds the adjoint of the inverse of $A$,
when $A$ is a convolution operator or the Radon transform.
Section~\ref{funcsec} proves that the $\delta$-net estimator
achieves the optimal rate of convergence in the ellipsoidal framework
and
it contains an oracle inequality for additive models.
Section~\ref{proofsec} contains the proofs of the main results.
The appendix contains calculations related to ellipsoids. 
%and an introduction for students to the setting of the article.

\paragraph{Notation}
We use the notation $\|\cdot \|$ to mean the Euclidean norm in ${\bf
R}^d$. The $L_2$ norm of a function $g:{\bf R}^{d} \to {\bf R}$ will
be denoted by $\|g\|_2$. The unit sphere in ${\bf R}^{d}$ is denoted
by ${\bf S}_{d-1} = \{ x \in {\bf R}^{d} : \|x\|=1\}$. The Lebesgue measure on
${\bf S}_{d-1}$ is denoted by $\mu$. We will make use of the formula $\mu({\bf S}_{d-1})
= 2\pi^{d/2}/\Gamma(d/2)$.
By $I_{R}$ %$R\subset {\bf R}^d$,
we denote the indicator function, i.e.\ $I_R(x)=1$ when $x\in R$ and
$I_R(x)=0$ otherwise. 
We write $a_n \asymp b_n$ to mean that $0<
\liminf_{n\to\infty} a_n/b_n \leq  \limsup_{n\to\infty} a_n/b_n <
\infty$, and $a_n  \succcurlyeq b_n$ means that $\liminf_{n\to
\infty} a_n/b_n >0$. The Fourier transform of a function $g\in
L_1({\bf R}^d)$ is defined by
$$
({F}g)(\omega)
=  \int_{{\bf R}^{d}} \exp\{ i x^T\omega\} g(x) \, dx ,
\qquad
\omega \in {\bf R}^{d} ,
$$
where $i$ is the imaginary unit. We use also the notation ${F}_1g$
when $g:{\bf R} \to {\bf R}$ is univariate. We have
$$
g(x)
=
(2\pi)^{-d} \int_{{\bf R}^{d}} \exp\{-i x^T\omega\}
(Fg)(\omega) \, d\omega ,
\qquad
x\in {\bf R}^{d} .
$$
By Parseval's theorem, we have for $f,g \in L_1({\bf R}^d) \cap
L_2({\bf R}^{d})$,

$$
\int_{{\bf R}^{d}} fg
= (2\pi)^{-d} \int_{{\bf R}^{d}} ({F}f)({F}g) .
$$
Convolution of $f$ and $g$ is denoted by$f*g(x)=\int_{{\bf R}^d}
f(x-y)g(y)\,dy$. We have that
\begin{equation}  \label{convo-prop}
{F}(f*g)  = ({F}f) ({F}g) .
\end{equation}
The probability measures of the Gaussian white noise
process $Y_n$ and of the i.i.d.~sequence
$(Y_1,\ldots ,Y_n)$ are denoted by $P_{Af}^{(n)}$.
%The expectation
%with respect to this measures is denoted by $E_f$.

\section{$\delta$-net minimizer}  \label{deltanetsec}

\paragraph{Definition of the estimator}

Let ${\cal F}$ be a set of densities or signal functions $f:{\bf
R}^d \to {\bf R}$. Let ${\cal F}_{\delta}$ be a finite $\delta$-net
of ${\cal F}$ in the $L_2$ metric, where $\delta>0$. That is, for
each $f\in {\cal F}$ there is a $\phi\in {\cal F}_{\delta}$ such
that $\| f-\phi\|_2 \leq \delta$. Define the estimator $\hat{f}$ by
$$
\hat{f} = \mbox{argmin}_{\phi \in {\cal F}_{\delta}} \gamma_n(\phi) ,
$$
where $\gamma_n(\phi)$ is defined in (\ref{gamman}).
Typically we would like to choose a $\delta$-net of minimal cardinality.
We assume that ${\cal F}$  is bounded in the $L_2$ metric,
\begin{equation}  \label{B2def}
\sup_{g\in {\cal F}} \| g \|_2 \leq B_2 ,
\end{equation}
where $0<B_2<\infty$.

\paragraph{An upper bound to MISE}

Theorem \ref{ora-ent} gives
a bound for the mean integrated squared error of the estimate.
We may identify the first term in the bound as a bias term and the second
term as a variance term.
The variance term depends on
the operator norm of $Q$ over the $\delta$-net  ${\cal F}_{\delta}$.
We define this operator norm as
\begin{equation}  \label{opnormQ}
\varrho(Q,{\cal F}_{\delta})       %T_{\delta}
=
\max_{\phi, \phi' \in {\cal F}_{\delta},\phi\neq\phi'}
\, \frac{\left\| Q(\phi-\phi') \right\|_2}{\| \phi - \phi' \|_2} \, ,
\qquad
\delta >0,
\end{equation}
where $Q$ is defined by (\ref{q-def}).
In the case of density estimation we need the additional assumption
that $\varrho(Q,{\cal F}_{\delta}) \geq1$ and that
$A{\cal F}$ and $Q{\cal F}$
are bounded in the $L_{\infty}$ metric:
\begin{equation} \label{b2-binf}
\varrho(Q,{\cal F}_{\delta}) \geq 1,
\qquad
\sup_{f\in {\cal F}} \|Af\|_{\infty} \leq B_{\infty},
\qquad
\sup_{f\in {\cal F}} \|Qf\|_{\infty} \leq B_{\infty}',
\end{equation}
where $0< B_{\infty},B_{\infty}' <\infty$.
%The assumption of the $L_2$-boundedness of ${\cal F}$ is made to guarantee the
%existence of a finite $\delta$-net.

\begin{theorem}  \label{ora-ent}
For the density estimation we assume that (\ref{b2-binf}) is satisfied.
We have that for $f\in {\cal F}$,
$$
E \left\| \hat{f} - f\right\|_2^2
\leq
C_1 \delta^2
+ C_2 \, \frac{\varrho^2(Q,{\cal F}_{\delta})
\cdot (\log_e (\#{\cal F}_{\delta})+1)}{n} \,,
$$
where
\begin{equation} \label{c1}
C_1
=
(1-2\xi)^{-1} (1+2\xi) ,
\end{equation}
\begin{equation} \label{c2}
C_2
=
(1-2\xi)^{-1} \xi C_{\tau} ,
\end{equation}
\begin{equation}  \label{ctau-2}
C_{\tau} >0 ,
\end{equation}
and $\xi$ is such that
\begin{equation}  \label{xi-choice}
\left\{
\begin{array}{ll}
C_{\tau}^{-1} \left( 4B_{\infty}'/3 +
\sqrt{2\left[8(B_{\infty}')^2/9+C_{\tau}B_{\infty}\right]} \right)
\leq
\xi < 1/2,
& \mbox{ density estimation }
\\
\sqrt{2/C_{\tau}} \leq  \xi < 1/2,
& \mbox{ white noise. }
\end{array} \right.
\end{equation}
\end{theorem}

A proof of Theorem~\ref{ora-ent} is given in Section \ref{ora-ent-proof}.

\begin{remark}{\em
Theorem~\ref{ora-ent} shows that the $\delta$-net estimator achieves the
rate of convergence $\psi_n$, when $\psi_n$ is the solution of the equation
\begin{equation}  \label{rateeq-univ}
\psi_n^2 \asymp n^{-1} \varrho^2(Q,{\cal F}_{{\psi}_n})  %T_{\psi_n}^2
\log (\# {\cal F}_{\psi_n}) .
\end{equation}
%In the optimal case the net ${\cal F}_{\delta}$
%is chosen so that its cardinality is minimal.
We calculate the rate under the assumptions that
$\log (\# {\cal F}_{\delta})$ and $\varrho(Q,{\cal F}_{\delta})$
increase polynomially as $\delta$ decreases:
we assume that one can find a $\delta$-net whose cardinality satisfies
$$
\log (\# {\cal F}_{\delta}) = C \delta^{-b}
$$
for some constants $b,C>0$
and we assume that
$$
\varrho(Q,{\cal F}_{\delta}) = C' \delta^{-a}
$$
for some $a,C'>0$ (in the direct case $a=0$ and $C'=1$).
Then (\ref{rateeq-univ}) can be written as
$
\psi_n^2 \asymp n^{-1} \psi_n^{-2a-b}
$
and
the rate of the $\delta$-net estimator is
\begin{equation}  \label{rate-nice}
\psi_n \asymp n^{-1/[2(a+1)+b]} .
\end{equation}
Let ${\cal F}$ be a set of $s$-smooth $d$-dimensional functions,
so that $b=d/s$.
Then the rate is
$$
\psi_n \asymp n^{-s/[2(a+1)s+d]},
$$
which gives for the direct case $a=0$ the classical rate
$\psi_n \asymp n^{-s/(2s+d)}$.
}\end{remark}

\section{A lower bound for MISE}  \label{lower-bound-section}

Theorem \ref{lowtheo} gives a lower bound for the mean integrated squared error
of any estimator, when estimating densities or signal functions
$f:{\bf R}^d \to {\bf R}$
in the function class ${\cal F}$.
Theorem \ref{lowtheo} holds also for nonlinear operators.

\begin{theorem}  \label{lowtheo}
Let $A$ be a possibly nonlinear operator.
Assume that for each sufficiently small $\delta>0$ we find
a finite set ${\cal D}_{\delta} \subset {\cal F}$ for which
\begin{equation}  \label{l2assu}
\min \{  \|f-g\|_2  :  f,g \in {\cal D}_{\delta}  ,
\,\,\,
f \neq g \}
\geq
C_0 \delta
\end{equation}
and
%$D_2$ is the $L_2$ distance and
\begin{equation}  \label{klassu}
\left\{
\begin{array}{ll}
\max \{  \|f-g\|_2  :  f,g \in {\cal D}_{\delta}  \}  \leq C_1 \delta,
& \mbox{ white noise, }
\\
\max \{  D_K(f,g)  :  f,g \in {\cal D}_{\delta}  \}  \leq C_1 \delta,
& \mbox{ density estimation, }
\end{array}  \right.
\end{equation}
where $D_K^2(f,g)=\int \log_e(f/g)\, f$ is the Kullback-Leibler distance,
%where $D_K^2(P,Q)=\int \log(dP/dQ)\,dP$ is the Kullback-Leibler distance,
%where $D_K$ is the Kullback-Leibler distance defined in (\ref{kldef}),
and $C_0$, $C_1$ are positive constants.
Denote
%$N_{\delta} = \#{\cal D}_{\delta}$ and
$$
\varrho_K(A,{\cal D}_{\delta})    %V_{\delta}
=
\left\{ \begin{array}{ll}
\frac{1}{\sqrt{2}} \,
\max_{f,g \in {\cal D}_{\delta},f\neq g} \frac{\|A(f-g)\|_2}{\|f-g\|_2},
& \mbox{ white noise, }
\\
\max_{f,g \in {\cal D}_{\delta},f\neq g} \frac{D_K(Af,Ag)}{\|f-g\|_2},
& \mbox{ density estimation. }
\end{array} \right.
$$
Let $\psi_n$ be such that
\begin{equation}  \label{minimaxeq}
\log_e (\#{\cal D}_{\psi_n})  % \gtrsim
\succcurlyeq   %> C_1
n \psi_n^2 \,  \varrho_K^2(A,{\cal D}_{\psi_n})  , %V_{\delta_n}^2
\end{equation}
where $a_n  \succcurlyeq b_n$ means that
$\liminf_{n\to \infty} a_n/b_n >0$.
%for sufficiently large $n$
Assume that
\begin{equation}  \label{epsi-infi}
\lim_{n\to \infty} n  \psi_n^2  \varrho_K^2(A,{\cal D}_{\psi_n})  = \infty  .
\end{equation}
%Assume that $\lim_{n\to \infty} N_{\delta_n} = \infty$.
Then,
$$
\liminf_{n \to \infty} \psi_n^{-2 }
\inf_{\hat{f}} \sup_{f\in{\cal F}} E \|f -\hat{f} \|_2^2
>0 ,
$$
where the infimum is taken over all estimators.
That is, $\psi_n$ is a lower bound for the minimax rate of convergence.
\end{theorem}

A proof of Theorem~\ref{lowtheo} is given in Section~\ref{lowtheo-proof}.

\begin{remark}{\em
Theorem \ref{lowtheo} shows that one can get a lower bound $\psi_n$ for
the rate of converge by solving the equation
\begin{equation}    \label{rateeq-univ-lower}
\psi_n^2   \,
\varrho_K^{2}(A,{\cal D}_{\psi_n})
\asymp
n^{-1}
\log_e (\#{\cal D}_{\psi_n} )  .
\end{equation}
The upper bound in Theorem~\ref{ora-ent} depends on the
operator norm of $Q$, defined in (\ref{opnormQ}),
whereas the lower bound depends on the operator norm of $A$.
Note also that the operator norm $\varrho(Q,{\cal F}_{\psi_n})$
is on the different side of the
equation in (\ref{rateeq-univ}) than the operator norm
$\varrho_K(A,{\cal D}_{\psi_n})$ in the
equation (\ref{rateeq-univ-lower}).
}\end{remark}

\begin{remark}{\em
In the density estimation case one can easily check assumptions
(\ref{klassu}) and (\ref{epsi-infi}) if one assumes that the
functions in $A{\cal D}_{\delta}$ are bounded and bounded away from
$0$. Then,
\begin{equation} \label{klbounds}
C' \cdot \|A(f-g)\|_2  \leq  D_K(Af,Ag)  \leq   C \cdot \|A(f-g)\|_2  .
\end{equation}
and (\ref{klassu}) and (\ref{epsi-infi}) follow by the corresponding
conditions with Hilbert norms instead of Kullback-Leibler distances.
}\end{remark}

\section{Dense  minimizer}  \label{deepsec}

The dense minimizer minimizes the empirical risk over the whole
function class ${\cal F}$. In contrast to the $\delta$-net estimator
the minimization is not restricted to a $\delta$-net. We call this
estimator ``dense minimizer'' because it is defined as a minimizer
over a possibly uncountable function class. The $\delta$-net
estimator is more widely applicable: it may be applied also to
estimate unsmooth functions and it may be applied when the operator
is severely ill-posed. The dense minimizer may be applied only for
relatively smooth cases (the entropy integral has to converge).
Because it works  without a restriction to a $\delta$-net we have
available a larger toolbox of numerical algorithms that can be
applied.

\paragraph{Definition of the estimator}

Let ${\cal F}$ be a collection of functions $f:{\bf R}^d \to {\bf R}$,
which are bounded in the $L_2$ metric as in (\ref{B2def}),
and
let the estimator $\hat{f}$
be a minimizer of the empirical risk over ${\cal F}$, up to $\epsilon>0$:
$$
\gamma_n(\hat{f}) \leq \mbox{inf}_{g \in {\cal F}} \gamma_n(g) + \epsilon ,
$$
where $\gamma_n(\phi)$ is defined in (\ref{gamman}).
For clarity, we present separate theorems for the Gaussian white noise model
and for the density estimation model.

\subsection{Gaussian white noise}

\paragraph{An upper bound to MISE}
Let ${\cal F}_{\delta}$, $\delta>0$, be a $\delta$-net of ${\cal F}$,
with respect to the $L_2$ norm.
Define
\begin{equation} \label{t-dense}
\varrho(Q,{\cal F}_{\delta})
=
\max\left\{
\frac{\| Q(f-g)\|_2}{\|f-g\|_2} \, :
f\in {\cal F}_{\delta},g\in {\cal F}_{2\delta},f\neq g \right\} ,
\qquad
\delta>0 ,
\end{equation}
where $Q$ is the adjoint of the inverse of $A$, defined by (\ref{q-def}).
%Denote with $N_{\delta}$ the cardinality of ${\cal F}_{\delta}$.
Define the entropy integral
\begin{equation} \label{ent-int}
G(\delta)
\stackrel{def}{=}
\int_0^{\delta} \varrho(Q,{\cal F}_{u})  \sqrt{\log_e (\#{\cal F}_u) } \, du ,
\qquad
\delta \in (0,B_2]  ,
\end{equation}
where $B_2$ is the $L_2$ bound defined by  (\ref{B2def}).

\begin{theorem}  \label{chain-white}
Assume that
\begin{enumerate}
\item
the entropy integral in (\ref{ent-int}) converges,
\item
$G(\delta)/\delta^2$
is decreasing on the interval $(0,B_2]$,
\item
$\varrho(Q,{\cal F}_{\delta}) = c\delta^{-a}$,
where $0 \leq a <1$ and $c>0$,
\item
$\lim_{\delta\to 0} G(\delta)\delta^{a-1} = \infty$,
%and $\varrho(Q,{\cal F}_{\delta}) \sqrt{\log_e (\#{\cal F}_{\delta})} \geq \sqrt{\log_e 2}$
%for sufficiently small $\delta$,
\item
$\delta \mapsto \varrho(Q,{\cal F}_{\delta}) \sqrt{\log_e(\#{\cal F}_{\delta})}$
is decreasing on $(0,B_2]$.
\end{enumerate}
Let $\psi_n$ be such that
\begin{equation}  \label{rateeq}
\psi_n^2 \geq C \, n^{-1/2} G(\psi_n) ,
\end{equation}
where $C$ is a positive constant,
and assume that $\lim_{n\to \infty} n\psi_n^{2(1+a)} = \infty$.
Then,  for $f\in {\cal F}$,
$$
E \left\| \hat{f} -f \right\|_2^2
\leq
C' \left( \psi_n^2 + \epsilon \right) ,
$$
for a positive constant $C'$,
for sufficiently large $n$.
\end{theorem}

A proof of Theorem~\ref{chain-white} is given in
Section~\ref{proof-chain-whitenoise}

\begin{remark}{\em
Assumption 5 is a technical assumption which is used to replace a
Riemann sum by an  entropy integral. We prefer to write the
assumptions in terms of the entropy integral in order to make them
more readable. }\end{remark}

\begin{remark}{\em
We may write $\varrho(Q,{\cal F}_{\delta})$ in a simpler way when
there exists minimal $\delta$-nets ${\cal F}_{\delta}$
which are nested:
$$
{\cal F}_{2\delta} \subset {\cal F}_{\delta} .
$$
Then we may define alternatively
$$
\varrho(Q,{\cal F}_{\delta})
=
\max_{f,g\in {\cal F}_{\delta},f\neq g}
\frac{\| Q(f-g)\|_2}{\|f-g\|_2} \, .
$$
}\end{remark}

\begin{remark}{\em      \label{dense-examples}
Theorem~ \ref{chain-white} and Theorem \ref{chain-density}
show that the rate of convergence
of the dense minimizer is the solution of the equation
\begin{equation}  \label{rateeq-dense}
\psi_n^2 = n^{-1/2} G(\psi_n) .
\end{equation}
To get the optimal rate the net ${\cal F}_{\delta}$
is chosen so that its cardinality is minimal.
In the polynomial case
one can find a $\delta$-net whose cardinality satisfies
$$
\log (\# {\cal F}_{\delta}) = C \delta^{-b}
$$
for some constants $b,C>0$
and the operator norm satisfies
$$
\varrho(Q,{\cal F}_{\delta}) = C' \delta^{-a}
$$
for some $a,C'>0$. (In the direct case $a=0$ and $C'=1$.)
Thus the entropy integral $G(\delta)$ is finite when
$\int_0^{\delta} u^{-a-b/2} \,du < \infty$,
which holds when
\begin{equation}  \label{entint-finite}
a+b/2 < 1 .
\end{equation}
Then (\ref{rateeq-dense}) leads to
$
\psi_n^2 \asymp n^{-1/2} \psi_n^{-a-b/2+1}
$
and the rate of the dense minimization estimator is
\begin{equation}  \label{rate-nice-dense}
\psi_n \asymp n^{-1/[2(a+1)+b]} .
\end{equation}
This is the same rate as the rate of the $\delta$-net estimator given
in (\ref{rate-nice}).
We have the following example.
Let ${\cal F}$ be a set of $s$-smooth $d$-dimensional functions,
so that $b=d/s$.
Then condition (\ref{entint-finite})
may be written as a condition for the smoothness
index $s$:
$$
s> \frac{d}{2(1-a)} \, .
$$
When the problem is direct, then $a=0$,
and we have the classical condition $s>d/2$.
The rate is
$
\psi_n \asymp n^{-s/[2(a+1)s+d]},
$
which gives for the direct case $a=0$ the classical rate
$\psi_n \asymp n^{-s/(2s+d)}$.
}\end{remark}

\subsection{Density estimation}

Let us call a $\delta$-bracketing net of ${\cal F}$
with respect to the $L_2$ norm
a set of pairs of functions
${\cal F}_{\delta} = \{ (g_j^L,g_j^U) : j=1,\ldots ,N_{\delta} \}$
such that
\begin{enumerate}
\item
$\| g_j^L-g_j^U\|_2 \leq \delta$, $j=1,\ldots ,N_{\delta}$,
\item
for each $g\in {\cal F}$ there is $j=j(g)\in \{1,\ldots ,N_{\delta}\}$
such that
$g_j^L \leq g \leq g_j^U$.
\end{enumerate}
%Let $N_{\delta}$ be the cardinality of a
%$\delta$-bracketing net ${\cal F}_{\delta}$ of ${\cal F}$.
%Let us write
%${\cal F}_{\delta}=\{ (g_j^L,g_j^U) : j=1,\ldots ,N_{\delta} \}$
Let us denote
${\cal F}_{\delta}^L=\{ g_j^L : j=1,\ldots ,N_{\delta} \}$
and
${\cal F}_{\delta}^U=\{ g_j^U: j=1,\ldots ,N_{\delta} \}$.
Define
\begin{equation}  \label{td-dens}
\varrho_{den}(Q,{\cal F}_{\delta})
=
\max\left\{ 
%T_{\delta}^{(1)} , T_{\delta}^{(2)} 
\varrho(Q,{\cal F}_{\delta}^L,{\cal F}_{\delta}^U) ,
\varrho(Q,{\cal F}_{\delta}^L,{\cal F}_{2\delta}^L)
\right\} ,
\end{equation}
where
$$
\varrho(Q,{\cal F}_{\delta}^L,{\cal F}_{\delta}^U)
=
\max\left\{
\frac{\| Q(g^U-g^L)\|_2}{\|g^U-g^L\|_2} \, :
g^L\in {\cal F}_{\delta}^L, \, g^U \in {\cal F}_{\delta}^U  \right\}
$$
and
$$
\varrho(Q,{\cal F}_{\delta}^L,{\cal F}_{2\delta}^L)
=
\max\left\{
\frac{\| Q(f-g)\|_2}{\|f-g\|_2} \, :
f\in {\cal F}_{\delta}^L,\, g\in {\cal F}_{2\delta}^L,\, f\neq g \right\} ,
$$
for $\delta>0$.
Define the entropy integral
\begin{equation} \label{ent-int-density}
G(\delta)
\stackrel{def}{=}
\int_0^{\delta} \varrho_{den}(Q,{\cal F}_u) 
\sqrt{\log_e (\# {\cal F}_u)} \, du ,
\qquad
\delta \in (0,B_2] ,
\end{equation}
where $B_2 = \sup_{f\in{\cal F}} \|f\|_2$.

\begin{theorem}  \label{chain-density}
We make the Assumptions 1-5 of Theorem~\ref{chain-white}
(with operator norm $\varrho_{den}(Q,{\cal F}_{\delta})$ in place of
$\varrho(Q,{\cal F}_{\delta})$),
and in addition we assume that
$\sup_{f\in {\cal F}} \|Af\|_{\infty} < \infty$,
$\sup_{g\in {\cal F}_{B_2}^L \cup {\cal F}_{B_2}^U} \|Qg\|_{\infty}
< \infty$,
and that the operator $Q$ preserves positivity
($g \geq 0$ implies that $Qg \geq 0$).
%where
%${\cal F}_{\delta}^L = \{ g^L : (g^L,g^U) \in {\cal F}_{\delta}\}$
%and
%${\cal F}_{\delta}^U = \{ g^U : (g^L,g^U) \in {\cal F}_{\delta}\}$.
Let $\psi_n$ be such that
\begin{equation}   \label{rateeq-density}
\psi_n^2 \geq C \, n^{-1/2} G(\psi_n) ,
\end{equation}
for a positive constant $C$,
and assume that $\lim_{n\to \infty} n\psi_n^{2(1+a)} = \infty$.
Then, for $f\in {\cal F}$,
$$
E \left\| \hat{f} -f \right\|_2^2
\leq
C' \left( \psi_n^2 + \epsilon \right) ,
$$
for a positive constant $C'$,
for sufficiently large $n$.
\end{theorem}

A proof of Theorem~\ref{chain-density} is given in
Section~\ref{proof-chain-density}. An analogous discussion of
optimal rates as in Remark~\ref{dense-examples} for the Gaussian
white noise model also applies for dense density estimators.

\section{Examples of operators}   \label{opesec}

As examples for operators we consider convolution operators and the
Radon transform. 
The definition of the empirical risk involves the adjoint of the
inverse of the operator $A$,
and we calculate the adjoint of the inverse of $A$, when $A$ is  a
convolution operator or the Radon transform.

%In Section \ref{opesec} 
%we calculate the adjoint of
%the inverse of these operators. 
%\subsection{Operators}  \label{opesec}

\subsection{Convolution}

The convolution operator $A$ is defined by
$$
Af = a*f ,
\qquad
f:{\bf R}^d \to {\bf R},
$$
where $a:{\bf R}^d\to {\bf R}$ is a known integrable function.
The adjoint of the inverse of $A$ is $Q$,
defined for $g:{\bf R}^d \to {\bf R}$, by
\begin{equation}  \label{convo-adjoint}
Qg  = F^{-1} \left( \frac{Fg}{Fa} \right)   ,
\end{equation}
where $F$ denotes the Fourier transform.
To derive this equation note that, for $h: {\bf R}^d \to {\bf R}$, 
$$
FA^{-1}h =  \frac{F h}{Fa}.
$$
Thus, for $h:{\bf R}^d \to {\bf R}$, $g:{\bf R}^d \to {\bf R}$,
applying  two times Parseval's theorem give
$$
\int_{{\bf R}^d}  (A^{-1}h) g
=
(2\pi)^d \int_{{\bf R}^d} \, \frac{(Fh) (Fg)}{Fa}
=
\int_{{\bf R}^d} h (Qg)  .
$$
Convolution operators appear in density estimation when the
observations contain additional measurement errors. In the
errors-in-variables model we observe $Y_i = X_i +\epsilon_i$,
$i=1,\ldots ,n$, where $X_i \sim f$, $f:{\bf R}^d \to {\bf R}$ is
the unknown density which we want to estimate, and $\epsilon_i \sim
a$ are the measurement errors. The density of the observations $Y_i$
is $Af = a * f$.

\subsection{Radon transform}

%\paragraph{Radon transform}
The Radon transform has been discussed in a series of papers and books including
\cite{Deans83}
and
\cite{Natte01}.
The Radon transform is defined as the integral of a $d$-dimensional
function over $d-1$-dimensional hyperplanes.
We parameterize the $d-1$-dimensional hyperplanes in
the $d$-dimensional Euclidean space
with the help of a direction vector
$\xi \in {\bf S}_{d-1}$
and a distance from the origin $u \in [0,\infty)$:
\begin{equation}  \label{plane-parameter}
P_{\xi,u} = \{ z \in {\bf R}^{d} : z^T\xi = u \} ,
\qquad
\xi \in {\bf S}_{d-1}, u \in [0,\infty)  .
\end{equation}
Define the Radon transform for $f:{\bf R}^{d} \to {\bf R}$ as
$$
(Af)(\xi,u) = \int_{P_{\xi,u}} f  ,
\qquad
\xi \in {\bf S}_{d-1}, \,\,\,  u \in [0,\infty)  ,
$$
where the integration is with respect to the $d-1$-dimensional
Lebesgue measure.
We will take the Radon transform as a mapping from
functions $f:{\bf R}^d \to {\bf R}$ to functions $Af : {\bf Y} \to {\bf R}$,
where ${\bf Y} = {\bf S}_{d-1} \times [0,\infty)$, and
the measure $\nu$ of the Borel space $({\bf Y}, {\cal Y}, \nu)$
is taken to be $d\nu(\xi,u) = u^{d-1} \, du \, d\mu(\xi)$.

The adjoint of the inverse of $A$ is $Q$,
defined for $g:{\bf R}^d \to {\bf R}$, by
\begin{equation}  \label{defperusidea}
(Qg)(\xi,u)  =
(2\pi)^{d-1}  \cdot (F_1^{-1} {\cal I}_{\xi} g) (u) ,
\qquad
\xi \in {\bf S}_{d-1},  \,\,\, u \in [0,\infty),
\end{equation}
where
$$
({\cal I}_{\xi} g)(t) = (Fg)(t\xi) ,
\qquad
\xi \in {\bf S}_{d-1},  \,\,\, t \in [0,\infty) .
$$
To see this note first that, 
for $h:{\bf S}_{d-1}\times [0,\infty) \to {\bf R}$,
we have that
\begin{equation}  \label{perusidea}
(FA^{-1}  h)(\omega)
=
({\cal H}_{\omega/\|\omega\|} h)(\|\omega\|) ,
\qquad
\omega \in {\bf R}^d,
\end{equation}
where ${\cal H}_{\xi}$ is the Fourier transform of $h(\xi,\,\cdot\,)$
for fixed $\xi\in {\bf S}_{d-1}$:
$$
{\cal H}_{\xi} h = F_1(h(\xi,\,\cdot\,)) ,
\qquad
\xi \in {\bf S}_{d-1} .
$$
Equation (\ref{perusidea}) follows directly
from the projection theorem, see \cite{Natte01}.

Two applications of Parseval's theorem and (\ref{perusidea}) 
give for $h:{\bf S}_{d-1} \times [0,\infty)  \to {\bf R}$,
$g:{\bf R}^d \to {\bf R}$, that
\begin{eqnarray*}
\int_{{\bf R}^d}  (A^{-1}h) g
& = & 
(2\pi)^d  \int_{{\bf R}^d}
( {\cal H}_{\omega/\|\omega\|} h)(\|\omega\|)
(Fg)(\omega) \,d\omega
\\ & = & 
 (2\pi)^d  \int_{{\bf S}_{d-1}} \int_0^{\infty} t^{d-1} ( {\cal H}_{\xi} h)(t)
(Fg)(t\xi) \,dt \, d\mu(\xi)
\\ & = &
 (2\pi)^{d-1}  \int_{{\bf S}_{d-1}} \int_0^{\infty} u^{d-1}
h(\xi ,u) (F_1^{-1} I_{\xi} g)(u) \,du \, d\mu(\xi)
\\ & = &
\int_{\bf Y} h (Qg) .
\end{eqnarray*}
This shows (\ref{defperusidea}).

%We applied  the projection theorem,
%which is also called the ``central slice theorem''
%or  the ``projection-slice theorem''.
%The projection theorem states that for $f: {\bf R}^d \to {\bf R}$,
%$g(u) = (Af)(\xi,u)$, $u\in {\bf R}$  ,
%\begin{equation}  \label{proje-theo}
%(F_1g)(t)
%=
%(Ff)(t\xi),
%\qquad
%t \in {\bf R}.
%\end{equation}

\paragraph{2D Radon transform}

In the 2D case we consider reconstructing a 2-di\-men\-si\-onal function
from observations of its integrals over lines.
Let $D = \{ x \in {\bf R}^2 : \|x\| \leq 1 \}$ be the unit disk
in ${\bf R}^2$.
The plane in (\ref{plane-parameter}) can be written as
$P_{\xi,u}=\{ u\xi + t\xi^{\perp} : t\in {\bf R}\}$,
where $\xi^{\perp}$ is a vector which is orthogonal to $\xi$.
We can write $\xi = (\cos\phi,\sin\phi)$ and
$\xi^{\perp} = (-\sin\phi,\cos\phi)$.
Thus we parameterize the lines by the length $u\in [0,1]$
of the perpendicular from the origin to the line
and by the orientation $\phi \in [0,2\pi)$ of this
perpendicular.
A common way to define 2D Radon transform is
\begin{equation}  \label{2d-radon}
Af(u,\phi)
=
\frac{\pi}{2\sqrt{1-u^2}} \,
\int_{\sqrt{1-u^2}}^{\sqrt{1-u^2}} f(u\cos\phi - t\sin\phi,
u\sin\phi+t\cos\phi) \, dt,
\end{equation}
where
$(u,\phi) \in {\bf Y} = [0,1]\times [0,2\pi]$,
and
we suppose that $f\in L_1(D) \cap L_2(D)$.
Now the Radon transform is $\pi$ times the average of $f$ over the line
segment that intersects $D$.
We consider $Rf$ as the element of $L_2({\bf Y},\nu)$,
where
%$S = [0,1]\times [0,2\pi]$ and
$\nu$ is the measure defined by
$d\nu(u,\phi) = 2\pi^{-1}\sqrt{1-u^2} \, du\, d\phi$.

\paragraph{Tomography}

The positron emission tomography is a density estimation
problem but the X-ray tomography is a regression type problem.
In the setting of positron emission tomography events happen at points
$X_1,\ldots ,X_n \in {\bf R}^d$,
and these points are i.i.d.~with density $f$.
We do not observe the location of the points but only that an event has
occurred on a hyperplane containing the point.
We assume that the hyperplane is uniformly oriented,
and that the distance of the hyperplane from the origin
is given by the Radon transform:
\begin{equation}  \label{u-cond-s}
S \sim \mbox{Unif}({\bf S}_{d-1}),
\qquad
U \, | \,  S=\xi
\sim
(Af) (\xi,\cdot) ,
\end{equation}
where hyperplanes are written as
$\{ z \in {\bf R}^d  : z^TS = U\}$.
We assume to observe i.i.d random variables
$Y_i = (S_i,U_i)\in {\bf S}_{d-1} \times [0,\infty)$, $i=1,\ldots ,n$,
which are distributed as $(S,U)$,
This is equivalent to observing the hyperplanes
$\{ z \in {\bf R}^d  : z^TS_i = U_i\}$.
We want to estimate the density
$f:{\bf R}^d \to {\bf R}$ in (\ref{u-cond-s}).
The density of the observations $Y_i$ is equal to
\begin{equation}  \label{tomo-ope}
(\tilde{A}f)(\xi,u) = \frac{1}{\mu({\bf S}_{d-1})} \, (Af)(\xi,u),
\qquad
\xi \in {\bf S}_{d-1}, \,\, u \in [0,\infty) .
\end{equation}
%We have
%$$
%{\bf Y} = {\bf S}_{d-1} \times [0,\infty) ,
%\qquad
%\| h \|_{2,\mu}^2
%= \int_{{\bf S}_{d-1}} \int_0^{\infty} u^{d-1} h^2(\xi,u) \,du\,d\mu(\xi),
%$$
%here $h :{\bf Y} \to {\bf R}$.

\section{Examples of function spaces}   \label{funcsec}

\subsection{Ellipsoidal function spaces}

%\subsection{Operator norms in the ellipsoidal framework}   \label{funcsec}

Since we
are in the $L_2$ setting it is natural to work in the sequence
space; we define the function classes as ellipsoids.
We shall apply singular value decompositions of the operators
and 
wavelet-vaguelette systems in the calculation of the rates
of convergence. In Section~\ref{ellipsoid-framework} we calculate
the operator norms in the framework of singular value decomposition. 
In Section~\ref{wave-vagu} we calculate the operator norms in the
wavelet-vaguelettte framework.
Section \ref{ratesec} derives the
rate of convergence of the $\delta$-net estimator for the case of a
convolution operator and the Radon transform, and the lower bound
for the rate of convergence of any estimator.

\subsubsection{Singular value decomposition}%{Ellipsoidal framework}
\label{ellipsoid-framework}

We assume that the underlying function space ${\cal F}$ consists of
$d$-variate functions that are linear combinations of orthonormal
basis functions $\phi_{j}$ with multi-index $j=(j_1,\ldots ,j_d) \in
\{0,1,\ldots\}^d$.
Define the ellipsoid
and
the corresponding collection of functions by
\begin{equation}  \label{sobolev-ellipsoid}
\Theta
=
\left\{ \theta : \sum_{j_1=0,\ldots ,j_d=0}^{\infty}
a_j^2 \theta_j^2 \leq L^2 \right\}  ,
\,\,
{\cal F}
=
\left\{ \sum_{j_1=0,\ldots ,j_d=0}^{\infty}  \theta_j\phi_j  :
\theta \in \Theta \right\}  .
\end{equation}

\paragraph{$\delta$-net and $\delta$-packing set for polynomial ellipsoids}
\label{sobolev-class}  \label{anisotropic}

We assume that there exists positive constants $C_1,C_2$ such that
for all $j \in \{0,1,\ldots\}^d$
\begin{equation}  \label{abound}
C_1\cdot  | j |^s \leq  a_j   \leq   C_2\cdot | j |^s ,
\end{equation}
where $| j | = j_1+\cdots + j_d$.
We construct a $\delta$-net
$\Theta_{\delta}$ and a $\delta$-packing set $\Theta_{\delta}^*$
in Appendix~\ref{ellipsoid-appendix}.
Since the construction is in the sequence space
we define the $\delta$-net and $\delta$-packing set of ${\cal F}$ by
\begin{equation}   \label{delta-net-definition-0}
{\cal F}_{\delta}
=
\left\{
\sum_{j_1=0,\ldots ,j_d=0}^{\infty}
\theta_j\phi_j  :  \theta \in \Theta_{\delta} \right\}  ,
\,
{\cal D}_{\delta}
=
\left\{
\sum_{j_1=0,\ldots ,j_d=0}^{\infty}
\theta_j\phi_j  :  \theta \in \Theta_{\delta}^* \right\}  .
\end{equation}
The set $\Theta_{\delta}$ is such that for
$\theta \in \Theta_{\delta}$
$$
\theta_j = 0,
\qquad
\mbox{ when } j \notin \{1,\ldots ,M \}^d,
$$
where
\begin{equation}   \label{M-definition}
M
\asymp
\delta^{-1/s} .
%\lceil ( C_1^{-1} 2^{1/2} L\delta^{-1})^{1/s}\rceil,
\end{equation}
%where $\lceil x\rceil$ denotes the smallest integer $\geq x$.
Set $\Theta_{\delta}^*$ is such that
for all $\theta \in \Theta_{\delta}^*$
\begin{equation}  \label{truncation}
\theta_j = \theta_j^* ,
\qquad
\mbox{ when } j \notin \{M^*,\ldots ,M \}^d,
\end{equation}
where
$\theta^*$ is a fixed sequence with
$\sum_{|j|\geq 0}^{\infty} a_j^2 {\theta_j^*}^2 = L^* < L$,
$$  %\begin{equation}   \label{M-definition}
M^* = [ M /2] .
$$ %\end{equation}
Furthermore, it holds that
\begin{equation}  \label{ellipsoid-cardinality}
\log( \# \Theta_{\delta})
\leq
C \delta^{-d/s} ,
\qquad
\log( \# \Theta_{\delta}^*)
\geq
C' \delta^{-d/s} .
\end{equation}

\paragraph{Operator norms}

We calculate the operator norms $\varrho(Q,{\cal F}_{\delta})$ and
$\varrho_K(A,{\cal D}_{\delta})$ in the ellipsoidal framework, where
${\cal F}_{\delta}$ and ${\cal D}_{\delta}$
are defined in (\ref{delta-net-definition-0}) and
Appendix~\ref{ellipsoid-appendix}.
We apply the singular value decomposition of $A$. We assume that  the
domain of $A$ is a separable Hilbert space $H$ with inner product
$\langle \cdot,\cdot \rangle$. 
The underlying function space ${\cal F}$ satisfies ${\cal F} \subset H$.
 We denote with $A^*$ the adjoint of $A$.
We assume that $A^*A$  is a compact operator on $H$ with eigenvalues $(b_j^2)$,
$b_j>0$, $j\in \{ 0,1,\ldots \}^d$, with orthonormal system of eigenfunctions
$\phi_j$.
We assume that there exists positive constants $q$ and $C_1,C_2$ such that
for all $j \in \{0,1,\ldots\}^d$
\begin{equation}   \label{bk-assumption}
C_1\cdot  | j |^{-q}  \leq  b_j   \leq   C_2\cdot | j |^{-q} .
\end{equation}
Let $g,g'$ in ${\cal F}_{\delta}$ or in ${\cal D}_{\delta}$,
respectively. Write
$$
g-g' = \sum_{j_1=1,\ldots ,j_d=1}^{\infty} (\theta_j-\theta_j')
\phi_j .
$$

\begin{enumerate}
\item
The functions $Q\phi_j$ are orthogonal and $\| Q\phi_j \|_2 = b_j^{-1}$.
Indeed, $Q = (A^{-1})^*$, and thus
$$
\langle Q\phi_j , Q\phi_l \rangle
=
\langle \phi_j , A^{-1}(A^{-1})^* \phi_l \rangle
=
b_l^{-2} \langle \phi_j , \phi_l\rangle  ,
$$
where we used the fact
\footnote{
Note that 
when a bounded linear operator $A$ between Banach spaces has a bounded inverse,
then $(A^{-1})^* = (A^*)^{-1}$, see
\cite{DunfoSchwa58}, Section~VI, Lemma~7, page~479.
}
$$
A^{-1} (A^{-1})^* \phi_l
=
A^{-1}(A^*)^{-1} \phi_l
=
(A^*A)^{-1} \phi_l
=
b_l^{-2} \phi_l .
$$
Thus for $g,g'\in  {\cal F}_{\delta}$,
\begin{eqnarray}
\| Q(g - g')\|_2^2
& = &      \nonumber
\left\|
\sum_{j_1=0,\ldots ,j_d=0}^{M} (\theta_j - \theta_j')^2 Q \phi_j \right\|_2^2
\\ & = &  \nonumber
\sum_{j_1=0,\ldots ,j_d=0}^{M}
(\theta_j - \theta_j')^2 b_j^{-2}
\\ & \leq &
C M^{2q}
\sum_{j_1=0,\ldots ,j_d=0}^{M}
(\theta_j - \theta_j')^2  ,
\label{ope-calcu}
\end{eqnarray}
where we used (\ref{bk-assumption}) to infer that
when $j \in \{ 0,\ldots ,M\}^d$, then
$$
b_j^{-2} \leq C_1^{-2} \cdot |j|^{2q} \leq C_1^{-2}  \cdot (dM)^{2q} .
$$
On the other hand,
$\|g-g'\|_2² = \sum_{j_1=0,\ldots ,j_d=0}^{M} (\theta_j-\theta_j')^2$.
This gives the upper bound for the operator norm
\begin{equation}   \label{Q-bound}
\varrho(Q,{\cal F}_{\delta})
\leq
C M^{q}
\leq
C' \delta^{-q/s}  ,
\end{equation}
by the definition of $M$ in (\ref{M-definition}).

\item
The functions  $A\phi_j$ are orthogonal and $\| A\phi_j \|_2 = b_j$.
Indeed,
$$
\langle A\phi_j , A\phi_l \rangle
=
\langle \phi_j , A^*A \phi_l \rangle
=
b_l^2 \langle  \phi_j , \phi_l \rangle  .
$$
Thus for $g,g'\in  {\cal D}_{\delta}$,
\begin{eqnarray*}
\| A(g - g')\|_2^2
& = &
\sum_{j_1=M^*,\ldots ,j_d=M^*}^{M}
(\theta_j - \theta_j')^2 \| A \phi_j\|_2^2
\\ & = &
\sum_{j_1=M^*,\ldots ,j_d=M^* }^{M} (\theta_j - \theta_j')^2 b_j^{2}  .
\end{eqnarray*}
This %combined with (\ref{truncation})
and similar calculations as in (\ref{ope-calcu})
imply that
\begin{equation}   \label{rho-lower}
C' \delta^{q/s} \leq \varrho_K(A,{\cal D}_{\delta}) \leq  C \delta^{q/s}  .
\end{equation}

\end{enumerate}

\subsubsection{Wavelet-vaguelette decomposition} \label{wave-vagu}

We assume that the underlying function space ${\cal F}$
consists of $d$-variate functions which are linear combinations
of orthonormal wavelet functions $(\phi_{jk})$,
where
$j \in \{0,1,\ldots \}$
and
$k \in \{0,\ldots ,2^j-1\}^d$.
%Define the Besov body,
%for
%$s = \sigma-(d/2-d/p)_+$,
%$\sigma > (d/p-d/2)_+$,
%$p \geq 1$,
%$0<q\leq \infty$,
%$L>0$,
%as
%$$
%\Theta
%=
%\left\{
%\begin{array}{ll}
%\left\{ \theta :
%\sum_j \left[ 2^{sj} \left( \sum_k |\theta_{jk}|^p \right)^{1/p} \right]^q
%\leq L^q
%\right\} ,
%\qquad  &
%q < \infty  ,
%\\
%\left\{ \theta :
%\sup_j 2^{sj} \left( \sum_k |\theta_{jk}|^p \right)^{1/p} \leq L
%\right\} ,
%\qquad   &
%q = \infty  .
%\end{array}
%\right.
%$$
%We shall restrict ourselves to the case $p=q=2$.
%In this case $\Theta$ is an $l_2$-body
The $l_2$-body and the corresponding class of functions
can now be defined as
$$
\Theta
=
\left\{ \theta :
\sum_j 2^{2sj} \sum_k |\theta_{jk}|^2
\leq L^2
\right\} ,
\qquad
{\cal F}
=
\left\{ \sum_{j}\sum_k  \theta_{jk}\phi_{jk}  :
\theta \in \Theta \right\}  ,
$$
where
$s>0$.
We have already constructed a $\delta$-net and $\delta$-packing
set for the $l_2$-bodies in (\ref{delta-net-definition-0}),
but in the current setting
for $\theta \in \Theta_{\delta}$
$$
\theta_{jk} = 0,
\qquad
\mbox{ when } j \geq J+1,
$$
where
\begin{equation}   \label{J-definition}
2^J
\asymp
\delta^{-1/s}
\end{equation}
and for $\theta \in \Theta_{\delta}^*$
$$
\theta_{jk} = \theta_{jk}^*,
\qquad
\mbox{ when } j \leq J^*
\mbox{ or }   j \geq J+1,
$$
where
$\theta^*$ is a fixed sequence with
$\sum_{j=0}^{\infty}\sum_k a_j^2 {\theta_{jk}^*}^2 = L^* < L$,
and
$ %\begin{equation}   \label{J-star-definition}
J^* =  J-1  .
$ %\end{equation}

\paragraph{Operator norms}

We can apply the wavelet-vaguelette decomposition,
as defined in \cite{Donoh95},
to calculate the operator norms
$\varrho(Q,{\cal F}_{\delta})$ and $\varrho_K(A,{\cal D}_{\delta})$.
We have available the following three sets of functions:
$(\phi_{jk})_{jk}$ is an orthogonal wavelet basis and
$(u_{jk})_{jk}$ and $(v_{jk})_{jk}$ are near-orthogonal sets:
$$
\left\| \sum_{jk} a_{jk} u_{jk} \right\|_2
\asymp
\| (a_{jk}) \|_{l_2} ,
\qquad
\left\| \sum_{jk} a_{jk} v_{jk} \right\|_2
\asymp
\| (a_{jk}) \|_{l_2} ,
$$
where $a \asymp b$ means that there exists positive constants $C,C'$ such that
$C b \leq a \leq C' b$.
The following quasi-singular relations hold:
$$
A\phi_{jk}
=
\kappa_j v_{jk},
\qquad
A^*u_{jk}
=
\kappa_j \phi_{jk} ,
$$
where $\kappa_j$ are quasi-singular values.
We assume that there exists positive constants $q$ and $C_1,C_2$ such that
for all $j \in \{0,1,\ldots\}$
\begin{equation}   \label{lambda-assumption}
C_1\cdot  2^{-qj}  \leq  \kappa_j   \leq   C_2\cdot 2^{-qj} .
\end{equation}

\begin{enumerate}
\item
Let $g,g' \in {\cal F}_{\delta}$.
Write
$$
g-g'
=
\sum_{j=0}^J \sum_{k} (\theta_{jk}-\theta_{jk}') \phi_{jk} .
$$
Since $Q = (A^{-1})^*$, then $Q A^* =  (A A^{-1})^* = I$.
Thus,  %denoting $\lambda=(j,k)$, $\lambda'=(j',k')$,
\begin{eqnarray*}
\langle Q\phi_{jk} , Q\phi_{j'k'} \rangle
& = &
\kappa_j^{-1} \kappa_{j'}^{-1} \langle Q A^* u_{jk} , Q A^* u_{j'k'} \rangle
\\ & = &
\kappa_j^{-1} \kappa_{j'}^{-1} \langle  u_{jk} , u_{j'k'} \rangle .
\end{eqnarray*}
Thus,
\begin{eqnarray}
\| Q(g - g')\|_2^2
& = &   \nonumber
\left\|
\sum_{j=0}^{J} \sum_{k} (\theta_{jk}-\theta_{jk}') Q\phi_{jk}
\right\|_2^2
\\ & = &   \nonumber
\left\|
\sum_{j=0}^{J} \kappa_j^{-1}
\sum_{k} (\theta_{jk}-\theta_{jk}') u_{jk}
\right\|_2^2
\\ & \asymp &   \nonumber
\sum_{j=0}^{J} \kappa_j^{-2}  \sum_k
(\theta_{jk} - \theta_{jk}')^2
\\ & \leq &   \label{ope-uppers}
C 2^{2qJ}
\sum_{j=0}^{J} \sum_k
(\theta_{jk} - \theta_{jk}')^2  ,
\end{eqnarray}
where we used (\ref{lambda-assumption}) to infer that
when $j \in \{ 0,\ldots ,J\}$, then
$$
\kappa_j^{-2} \leq C_1^{-2} \cdot 2^{2qj} \leq C_1^{-2}  \cdot 2^{2qJ} .
$$
On the other hand,
$\|g-g'\|_2^2 = \sum_{j=0}^{J} \sum_k (\theta_{jk}-\theta_{jk}')^2$.
This gives the upper bound for the operator norm
$$ %\begin{equation}   \label{Q-bound}
\varrho(Q,{\cal F}_{\delta})
\leq
C 2^{qJ}
\leq
C' \delta^{-q/s}  ,
$$  %\end{equation}
by the definition of $J$ in (\ref{J-definition}).

\item
We have
$
\langle A\phi_{jk} , A\phi_{j'k'} \rangle
=
\kappa_j \kappa_{j'} \langle  v_{jk} , v_{j'k'} \rangle
$
and $(v_{jk})$ is a near-orthogonal set.
Thus, similarly as in (\ref{ope-uppers}), we get
$$  %\begin{equation}  \label{rho-lower}
C' \delta^{q/s}
\leq
\varrho_K(A,{\cal D}_{\delta})
\leq
C \delta^{q/s}  .
$$ %\end{equation}

\end{enumerate}
%\hspace*{\fill}  $\Box$

\subsubsection{Rates of convergence}    \label{ratesec}

%In Section  %\ref{sec-convolu} and in Section~\ref{sec-rado}
We derive the rates of convergence for the $\delta$-net estimator
when the operator is a convolution operator and the Radon transform.
It is also shown that the lower bounds have the same order as the upper bounds.
We give examples in the setting of the Gaussian white noise model.

\paragraph{Convolution}   %\label{sec-convolu}

Let $A$ be a convolution operator: $Af = a*f$
where $a :{\bf R}^d \to {\bf R}$ is a known function.
Denote
$$
\phi_{jk}(x)
=
\prod_{i=1}^d \sqrt{2} \left[ (1-k_i) \cos(2\pi j_ix_i) +
 k_i\sin(2\pi j_ix_i)  \right]  ,
\qquad
x\in [0,1]^d ,
$$
where
$j \in \{0,1,\ldots \}^d$, $k\in K_j$, where
$$
K_j
=
\left\{ k \in \{ 0,1\}^d  :  k_i=0, \mbox{ when } j_i=0   \right\}  .
$$
The cardinality of $K_j$ is $2^{d-\alpha(j)}$, where
$\alpha(j)= \#\{ j_i  :  j_i=0 \}$.
The collection $(\phi_{jk})$, $(j,k)\in \{0,1,\ldots \}^d \times K_j$,
is a basis for $1$-periodic functions on $L_2([0,1]^d)$.
When the convolution kernel $a$ is an $1$-periodic function in $L_2([0,1]^d)$,
then we can write
$$
a(x) = \sum_{j_1=0,\ldots ,j_d=0}^{\infty}  \sum_{k\in K_j}
b_{jk} \phi_{jk}(x) .
$$
The functions $\phi_{jk}$ are the singular
functions of the operator $A$ and the
values $b_{jk}$ are the corresponding singular values.
We assume that the underlying function space is equal to
\begin{equation}  \label{calf-definition}
{\cal F}
=
\left\{
\sum_{j_1=0,\ldots ,j_d=0}^{\infty}  \sum_{k\in K_j}
\theta_{jk} \phi_{jk}(x)
:
(\theta_{jk})   \in \Theta \right\}  ,
\end{equation}
where
\begin{equation}  \label{theta-ellipsoid}
\Theta
=
\left\{   \theta  :  \sum_{j_1=0,\ldots ,j_d=0}^{\infty}  \sum_{k\in K_j}
a_{jk}^2 \theta_{jk}^2 \leq L^2
\right\}  .
\end{equation}
We give the rate of convergence of the
$\delta$-net estimator and show that the estimator achieves the
optimal rate of convergence.
Optimal rates of convergence has been previously obtained for the
convolution problem in various settings in
\cite{Ermak89},
\cite{DonohLow92},
\cite{Koo93},
\cite{KorosTsyba93book}.
%\cite{Donoh95}.
%\cite{Johns99}.

\begin{corollary}  \label{coro-convo}
Let ${\cal F}$ be the function class as defined in (\ref{calf-definition}).
We assume that the coefficients of the ellipsoid (\ref{theta-ellipsoid})
satisfy
$$
C_0 |j|^s  \leq a_{jk} \leq C_1 |j|^s  .
$$
for some $s>0$ and $C_0,C_1>0$.
We assume that the convolution filter $a$ is 1-periodic function
in $L_2([0,1]^d)$
and that the Fourier coefficients of filter $a$ satisfy
$$
C_2 |j|^{-q} \leq b_{jk} \leq C_3 |j|^{-q}
$$
for some $q \geq 0$, $C_2,C_3>0$.
%Assume that $a$ satisfies (\ref{poly-dec}) with $s > q$.
Then,
$$
\limsup_{n\to\infty} n^{2s/(2s+2q+d)} \sup_{f\in {\cal F}}
E_f \left\| \hat{f} - f\right\|_2^2
< \infty ,
$$
where $\hat{f}$ is the $\delta$-net estimator.
%Assume that $a$ satisfies  (\ref{poly-dec}) and (\ref{poly-dec-lower}).
Also,
$$
\liminf_{n\to\infty} n^{2s/(2s+2q+d)} \inf_{\hat{g}} \sup_{f\in {\cal F}}
E_f \left\| \hat{g} - f\right\|_2^2
> 0 ,
$$
where the infimum is taken over any estimators $\hat{g}$.
\end{corollary}

{\em Proof.}
For the upper bound we apply Theorem \ref{ora-ent}.
Let ${\cal F}_{\delta}$ be the $\delta$-net of ${\cal F}$
as constructed in (\ref{delta-net-definition-0}).
We have shown in (\ref{Q-bound}) that
$$
\varrho(Q,{\cal F}_{\delta}) \leq C \delta^{-a} ,
$$
where $a = q/s$.
We have stated in (\ref{ellipsoid-cardinality}) that
the cardinality of the $\delta$-net satisfies
$$
\log (\# {\cal F}_{\delta})  \leq  C\delta^{-b},
$$
where $b=d/s$.
Thus we may apply (\ref{rate-nice}) to get the rate
$$
\psi_n =n^{-1/(2(a+1)+b)} = n^{-s/(2s+2q+d)}  .
$$
The upper bound is proved.
For the lower bound we apply Theorem \ref{lowtheo}.
Assumption (\ref{l2assu})
holds because ${\cal D}_{\delta}$ in
(\ref{delta-net-definition-0})
is a $\delta$-packing set.
Assumption (\ref{klassu})
holds by the construction, see (\ref{klassu-holds})
in Appendix~\ref{ellipsoid-appendix}.
Assumptions (\ref{minimaxeq})
and (\ref{epsi-infi}) follow from
(\ref{ellipsoid-cardinality})
and
(\ref{rho-lower}).
Thus the lower bound is proved.
\hspace*{\fill} $\Box$

\paragraph{Radon transform}  %\label{sec-rado}

We consider the 2D Radon transform as defined in (\ref{2d-radon}).
The singular value decomposition of the Radon transform can be found in
\cite{Deans83}.
Let
$$
\tilde{\phi}_{jk}(r,\theta)
=
\pi^{-1/2}(j+k+1)^{1/2} Z_{j+k}^{|j-k|}(r) e^{i(j-k)\theta},
\,\,\,
(r,\theta) \in D = [0,1] \times [0,2\pi) ,
$$
where
$Z_a^b$ denotes the Zernike polynomial of degree $a$ and order $b$.
Functions $\tilde{\phi}_{jk}$, $j,k=0,1,\ldots$, $(j,k)\neq (0,0)$,
constitute an orthonormal complex-valued basis for
$L_2(D)$.
The corresponding orthonormal functions in $L_2({\bf Y},\nu)$ are
$$
\tilde{\psi}_{jk}(u,\phi)
=
\pi^{-1/2} U_{j+k}(u) e^{i(j-k)\phi} ,
\qquad
(u,\phi) \in {\bf Y} = [0,1] \times [0,2\pi) ,
$$
where
$U_m(\cos \theta) = \sin((m+1)\theta)/\sin\theta$ are the
Chebyshev polynomials of the second kind.
We have
$$
A \tilde{\phi}_{jk} = b_{jk} \tilde{\psi}_{jk} ,
$$
where the singular values are
\begin{equation}  \label{radon-singular-values}
b_{jk}
=
\pi^{-1} (j+k+1)^{-1/2}  .
\end{equation}
We shall identify the complex bases with the equivalent real
orthonormal bases by
$$
\phi_{jk}
 =
\left\{   \begin{array}{ll}
\sqrt{2}\, \mbox{Re}(\tilde{\phi}_{jk})  &  \mbox{ if } j>k  \\
\tilde{\phi}_{jk}  &  \mbox{ if } j=k  \\
\sqrt{2}\,  \mbox{Im}(\tilde{\phi}_{jk})  &  \mbox{ if } j<k  .
\end{array}  \right.
$$
We assume that the underlying function space is equal to
\begin{equation}  \label{calf-definition-radon}
{\cal F}
=
\left\{
\sum_{j_1=0,j_2=0,(j_1,j_2)\neq(0,0)}^{\infty}
\theta_{j_1j_2} \phi_{j_1j_2}(x)
:
(\theta_{j_1j_2})   \in \Theta \right\}  ,
\end{equation}
where
\begin{equation}  \label{theta-ellipsoid-radon}
\Theta
=
\left\{   \theta  :  \sum_{j_1=0,j_2=0,(j_1,j_2)\neq (0,0)}^{\infty}
a_{j_1j_2}^2 \theta_{j_1j_2}^2 \leq L^2
\right\}  .
\end{equation}
%where
%$$
%a_{j_1j_2}
%=
%(j_1+j_2+1)^s  .
%%(j_1+1)^{s} (j_2+1)^{s}  .
%$$
We give the rate of convergence of the
$\delta$-net estimator and show that the estimator achieves the
optimal rate of convergence.
Optimal rates of convergence have been previously obtained in
\cite{JohnsSilve90},
%\cite{KorosTsyba89},
\cite{KorosTsyba91},
\cite{DonohLow92},
\cite{KorosTsyba93book}.

\begin{corollary}
Let ${\cal F}$ be the function class as defined in
(\ref{calf-definition-radon}).
We assume that the coefficients of the ellipsoid (\ref{theta-ellipsoid-radon})
satisfy
$$
C_0 |j|^s  \leq a_{jk} \leq C_1 |j|^s  .
$$
for some $s>0$ and $C_0,C_1>0$.
Then, for $d=2$,
$$
\limsup_{n\to\infty} n^{2s/(2s+2d-1)} \sup_{f\in {\cal F}}
E_f \left\| \hat{f} - f\right\|_2^2
< \infty .
$$
where $\hat{f}$ is the $\delta$-net estimator.
%Assume that $a$ satisfies  (\ref{poly-dec}) and (\ref{poly-dec-lower}).
Also,
$$
\liminf_{n\to\infty} n^{2s/(2s+2d-1)} \inf_{\hat{g}} \sup_{f\in {\cal F}}
E_f \left\| \hat{g} - f\right\|_2^2
> 0 ,
$$
where the infimum is taken over any estimators $\hat{g}$.
\end{corollary}

{\em Proof.}
For the upper bound we apply Theorem \ref{ora-ent}.
Let ${\cal F}_{\delta}$ be the $\delta$-net of ${\cal F}$
as constructed in (\ref{delta-net-definition-0}).
We have shown in (\ref{Q-bound}) that
$$
\varrho(Q,{\cal F}_{\delta}) \leq C \delta^{-a} ,
$$
where $a=q/s$ and $q=1/2$ (so that $a=(d-1)/(2s)$),
since the singular values are given in (\ref{radon-singular-values}).
We have stated in (\ref{ellipsoid-cardinality}) that
the cardinality of the $\delta$-net satisfies
$$
\log (\# {\cal F}_{\delta})  \leq  C\delta^{-b},
$$
where $b=d/s$.
Thus we may apply (\ref{rate-nice}) to get the rate
$$
\psi_n = n^{-s/(2s+2d-1)}  .
$$
The upper bound is proved.
For the lower bound we apply Theorem \ref{lowtheo}
similarly as in the proof of Corollary~\ref{coro-convo}.
\hspace*{\fill} $\Box$

\subsection{Additive models}

In this section we will show that our approach can be used to prove 
oracle results for additive models.  
In additive models the unknown 
function $f:{\bf R}^d \to {\bf R}$ is assumed to have an additive 
decomposition $f(x)=f_1(x_1)+\cdots+f_d(x_d)$ with unknown additive 
components $f_j: {\bf R} \to {\bf R}$, $j=1,\ldots,d$. We compare this 
model with theoretical oracle models where only one component function 
$f_r$ is unknown, but the other functions $f_j$ ($j \not = r$) are known. 
We will show below that the function $f$ can be estimated with the same 
rate of convergence as in the oracle model that has the slowest rate of 
convergence. In particular, if the rate of convergence is the same 
in all oracle models then the rate in the additive model remains the same. 
This is a well known fact for classical additive regression models, 
see e.g.~\cite{Stone85}. %Stone (1985). 
It efficiently avoids the curse
of dimensionality in contrast to the full dimensional
nonparametric model. Furthermore, it is practically important because
it allows a flexible and nicely interpretable model for regression
with high dimensional covariates, 
see e.g.~\cite{HastiTibsh90} % Hastie and Tibshirani (1990) 
for a discussion of
the additive and related models. 
Thus, our result will generalize the oracle result for additive models of 
\cite{Stone85} to inverse problems. 
For a theoretical discussion we will first use a slightly more general 
framework. 
We will come back to additive models afterwards.

\subsubsection{Abstract setting}

We assume that the function class $\cal F$ is a subset of the direct sum 
of spaces ${\cal F}_1$,\ldots, ${\cal F}_p$. 
All spaces contain functions from 
$f:{\bf R}^d \to {\bf R}$.
 At this stage, we do not assume that functions in 
${\cal F}_j$ ($j=1,\ldots,p$) depend only on the argument $x_j$. 
An example of this more general set up are sums of smooth functions and 
indicator functions of convex sets or of sets with smooth boundary. 
We assume that a finite $\delta$-net ${\cal F}_{\delta}$ of ${\cal F}$ 
is a subset of the direct sum 
${\cal F}_{1,\delta} \oplus \cdots \oplus {\cal F}_{p,\delta}$, 
where ${\cal F}_{j,\delta}$ are finite subsets of ${\cal F}_{j}$. 
We denote the number of elements of ${\cal F}_{j,\delta}$ by 
$\exp(\lambda_j)$. 
Furthermore, we write $\rho_j= \rho(Q, {\cal F}_{j,\delta})$. 
We make the following essential geometrical assumption:
\begin{equation} \label{geo1}
\|f_1+\cdots+f_p\|_2^2 \geq c \sum_{j=1}^p \|f_j\|_2^2 \end{equation}
for a positive constant $c > 0$.
For the $\delta$-net minimizer $\hat f$ over the $\delta$-net  
${\cal F}_{\delta}$ we get the following result in the white noise model. 
(An additive model for density estimation would not make much sense.)

\begin{theorem}   \label{thm4}  \label{ora-add}
We make assumption (\ref{geo1}). 
In the white noise model the following bound holds for the $\delta$-net 
minimizer $\hat f$, for $f \in {\cal F}$,
\begin{eqnarray*}
 E \left(\|\hat f - f\|^2_2 \right) 
\leq 
3 \delta^2 
+ 32 c^{-1} n ^{-1} \left[ \sum_{j=1}^p \rho_j^2 \lambda_j 
+ \left(\sum_{j=1}^p \rho_j\right)^2 \right].
\end{eqnarray*}
\end{theorem}

A proof of Theorem~\ref{ora-add} is given in Section~\ref{ora-add-proof}.

\subsubsection{Application to additive models}

We now apply Theorem~\ref{thm4} for discussing additive models 
$f(x)=f_1(x_1)+\cdots+f_d(x_d)$. In $L_2({\bf R}^d)$ we have
$\|f_1+\cdots+f_d\|_2^2 = \sum_{j=1}^d \|f_j\|_2^2$, 
if the functions $f_j$ are normed such that $\int f_j(x_j)\, dx_j=0$. 
Thus (\ref{geo1}) holds trivially. 
Assumption (\ref{geo1}) also holds in other $L_2$-spaces with 
dominating measure differing from the Lebesgue measure. 
A discussion of condition (\ref{geo1}) for these classes can be found 
e.g.~\cite{MammeLintoNiels99}. % in Mammen, Linton and Nielsen (1999). 
See also \cite{BickeKlaasRitovWelln93}.
%Bickel, Klaassen,  Ritov and Wellner (1993).
Such $L_2$-spaces naturally arise in additive regression models. 
For a white noise model they come up if one assumes an additive model 
for transformed covariables. We assume that for the models  ${\cal F}_{j}$ 
one can find $\delta_j$-nets  ${\cal F}_{j,\delta_j}$ such that choosing 
$\delta_j = \psi_{n,j}$ with 
$$
\psi_{n,j}^2 \asymp n ^{-1} \rho^2(Q,  {\cal F}_{j,\psi_{n,j}}) 
\log (\# {\cal F}_{j,\psi_{n,j}})
$$
gives a rate optimal $\delta$-net minimizer in the model ${\cal F}_{j}$. Now,
${\cal F}_{\delta}={\cal F}_{1,\delta_1} \oplus 
\cdots \oplus {\cal F}_{d,\delta_d}$ 
is a $\delta$-net of ${\cal F}$ with $\delta=\sum_{j=1}^d \delta_j$.
From Theorem~\ref{thm4} we get that the $\delta$-net minimizer $\hat f$ 
over the net ${\cal F}_{\delta}$ achieves the rate $O(\psi_n)$ 
with $\psi_n= \max_{1\leq j \leq d} \psi_{n,j}$. This is just the 
type of result we called oracle result at the beginning of this section.

In general, the oracle result does not follow from Theorem~\ref{ora-ent}. 
The application of Theorem~\ref{ora-ent} leads to an assumption of the type
$$
n ^{-1} \max_{1 \leq j \leq d}  
\rho^2(Q,  {\cal F}_{j,\psi_{n,j}}) \times \max_{1 \leq j \leq d}
 \log (\# {\cal F}_{j,\psi_{n,j}})
=
O\left(\psi_n^2\right)
$$ 
whereas Theorem~\ref{thm4} only requires that 
$$
n ^{-1} \max_{1 \leq j \leq d} \left[ \rho^2(Q,  {\cal F}_{j,\psi_{n,j}}) 
\log (\# {\cal F}_{j,\psi_{n,j}}) \right]
=
O\left(\psi_n^2\right).
$$
This can make a big difference. 
First of all the entropy numbers of the additive classes  ${\cal F}_{j}$ 
may differ. Furthermore, the operator $Q$ may act quite differently 
on the spaces ${\cal F}_{j}$.

\subsubsection{Ellipsoidal spaces and convolution}

As an example
we now assume that the underlying function space is
${\cal F} = {\cal F}_1 \oplus \cdots \oplus {\cal F}_d$,
where
$$
{\cal F}_k
=
\left\{
\sum_{j=0}^{\infty}
\theta_{kj}\phi_{kj}  :  \theta_{k\cdot} \in \Theta_{s_k,L_k} \right\}  
$$
for 
basis functions $\phi_{kj}:[0,1] \to {\bf R}$ 
and the ellipsoids are defined by
\begin{equation}  \label{sobolev-ellipsoid-additive-1}
\Theta_{s_k,L_k}
=
\left\{ \theta_{k\cdot} : \sum_{j=0}^{\infty}
a_{kj}^2 \theta_{kj}^2 \leq L_k^2 \right\}  ,
\qquad
k=1,\ldots ,d,
\end{equation}
where we assume that there exists positive constants $C_1,C_2$ such that
for all $j \in \{0,1,\ldots\}$
\begin{equation}  \label{abound-additive}
C_1\cdot  j^{s_k} \leq  a_{kj}   \leq   C_2\cdot j^{s_k} .
\end{equation}
Let $A$ be a convolution operator: $Af = a*f$
where $a :{\bf R}^d \to {\bf R}$ is a known function.
Then
$$
Af 
=
A_1f_1 + \cdots + A_df_d ,
$$
where
$f(x) = f_1(x_1) + \cdots + f_d(x_d)$
and
$$
A_kf_k(x_k)
=
\int_{[0,1]^d} f_k(x_k-y_k)a_k(y_k) \, dy_k ,
$$
where
$$
a_k(y_k)
=
\int_{[0,1]^d} a(y) \prod_{l=1,l\neq k}^d dy_l
$$
is the $k$th marginal function of $a$.
We can decompose $Q$ accordingly:
$$
Qg
=
Q_1g_1 + \cdots + Q_dg_d .
$$ 
Operators $A_j$ and $Q_j$ are restrictions of $A$ and $Q$
to ${\cal F}_j$.
We apply the singular value decomposition for $A_k$. 
Denote
$$
\phi_{kj}(t)
=
\sqrt{2} \cos(2\pi j t)   ,
\qquad
t\in [0,1] ,
$$
where
$j = 1,2,\ldots$ and $\phi_0(t) = I_{[0,1]}(t)$.
The collection $(\phi_{kj})$, $j = 0,1,\ldots$,
is a basis for $1$-periodic functions on $L_2([0,1])$.
When $a_k$ are $1$-periodic functions in $L_2([0,1])$,
then we can write
$$
a_k(x_k) = \sum_{j=0}^{\infty} b_{kj} \phi_{kj}(x_k) .
$$
The functions $\phi_{kj}$ are the singular
functions of the operator $A_k$ and the
values $b_{kj}$ are the corresponding singular values.
We give the rate of convergence of the
$\delta$-net estimator and show that the estimator achieves the
optimal rate of convergence.

\begin{corollary}  \label{coro-convo-additive}
Let ${\cal F} = {\cal F}_1 \oplus \cdots \oplus {\cal F}_d$.
We assume that the coefficients of the ellipsoid 
satisfy (\ref{abound-additive}).
We assume that $a_k$ are 1-periodic functions
in $L_2([0,1])$
and that the Fourier coefficients of $a_k$ satisfy
$$
C_2 j^{-q_k} \leq b_{kj} \leq C_3 j^{-q_k}
$$
for some $q_k \geq 0$, $C_2,C_3>0$.
Then, in the white noise model,
$$
\limsup_{n\to\infty} n^{a} \sup_{f\in {\cal F}}
E_f \left\| \hat{f} - f\right\|_2^2
< \infty ,
$$
where $\hat{f}$ is the $\delta$-net estimator
and 
$$
a = \min_{k=1,\ldots ,d} \frac{2s_k}{2s_k+2q_k+1} \,   .
$$
Also,
$$
\liminf_{n\to\infty} n^a \inf_{\hat{g}} \sup_{f\in {\cal F}}
E_f \left\| \hat{g} - f\right\|_2^2
> 0 ,
$$
where the infimum is taken over any estimators $\hat{g}$ 
in the white noise model.
\end{corollary}

{\em Proof.}
For the upper bound we apply Theorem~\ref{ora-add}.
As in Section~\ref{ellipsoid-framework} we can find 
$\delta$-nets ${\cal F}_{k,\delta}$ for ${\cal F}_k$ whose
cardinality is bounded by
$
\log ( \# {\cal F}_{k,\delta} )
\leq
C \delta^{-1/s_k} 
$
and
$
\varrho(Q_k,{\cal F}_{k,\delta}) \leq C \delta^{-q_k/s_k} .
$
The upper bound of Theorem~\ref{ora-add} gives as the rate
the maximum of the component rates 
$n^{-2s_k/(2s_k+2q_k+1)}$.
For the lower bound we 
apply the lower bound of Corollary~\ref{coro-convo} in the case $d=1$
and the fact that one cannot do better in the additive model
than in the model that has only one component.
\hspace*{\fill} $\Box$

\section{Proofs}  \label{proofsec}

\subsection{A preliminary lemma}

We prove that the theoretical error of a minimization estimator
may be bounded by the optimal theoretical error
and an additional stochastic term.

\begin{lemma}  \label{oralemma}
Let ${\cal C} \subset L_2({\bf R}^d)$.
Let $\hat{f}\in {\cal C}$
be such that
\begin{equation}  \label{glob}
\gamma_n(\hat{f}) \leq \inf_{g \in {\cal C}} \gamma_n(g)
+ \varepsilon,
\end{equation}
where $\varepsilon \geq 0$.
Then for each $f^0 \in {\cal C}$,
$$
\left\| \hat{f} -f \right\|_2^2
\leq
\left\| f^0 -f \right\|_2^2
+
\varepsilon
+
2 \nu_n[Q(\hat{f}-f^0)]
$$
where
$f$ is the true density or the true signal function,
and $\nu_n(g)$ is the centered empirical operator:
\begin{equation}  \label{nun-def}
\nu_n(g)
= \left\{  \begin{array}{ll}
\int g \, dY_n  - \int_{{\bf Y}} g (Af) ,
& \mbox{\em  white noise model, }
\\
n^{-1} \sum_{i=1}^n g(Y_i) - \int_{{\bf Y}} g (Af) ,
& \mbox{\em density estimation, }
\end{array} \right.
\end{equation}
where $g :{\bf R}^d \to {\bf R}$.
\end{lemma}

{\em Proof.}
We have for $g=\hat{f}$, $g=f^0$,
\begin{eqnarray*}
\lefteqn{
\| g-f\|_2^2 - \gamma_n(g)
} \\ & = &
\left\{  \begin{array}{ll}
\|f\|_2^2 - 2\int_{{\bf R}^d}fg  + 2 \int (Qg) \, dY_n ,
& \mbox{ white noise model }
\\
\|f\|_2^2 - 2\int_{{\bf R}^d}fg  + 2n^{-1} \sum_{i=1}^n(Qg)(Y_i) ,
& \mbox{ density estimation. }
\end{array} \right.
\end{eqnarray*}
We have $\int_{{\bf R}^d}fg = \int_{{\bf Y}}(Af)(Qg)$.
Thus,
\begin{equation}  \label{l2emp}
\left\| \hat{f}-f \right\|_2^2  - \gamma_n\left(\hat{f}\right)
+ \gamma_n\left(f^0\right) - \left\| f^0-f\right\|_2^2
=
2\nu_n\left[ Q\left(\hat{f} - f^0\right)\right] .
\end{equation}
Thus,
%since by the definition of $\tilde{f}_n$ as a minimization estimator,
%$\| \tilde{W}_n\|_{l_1} = \| \tilde{W}^0\|_{l_1} =1$,
\begin{eqnarray}
\lefteqn{
\left\| \hat{f} - f \right\|_2^2 - \left\| {f}^0 - f \right\|_2^2
} \nonumber \\ & = &  \nonumber
\left\| \hat{f} - f \right\|_2^2 - \gamma_n( \hat{f})
+ \gamma_n( \hat{f}) - \left\| {f}^0 - f \right\|_2^2
\\ & \leq & \label{do}
\left\| \hat{f} - f \right\|_2^2 - \gamma_n( \hat{f})
+ \gamma_n( {f}^0) +  \varepsilon
- \left\| {f}^0 - f \right\|_2^2
\\ & = &   \label{re}
2\nu_n\left[ Q\left(\hat{f} - f^0\right)\right]  +  \varepsilon  .
\end{eqnarray}
In (\ref{do}) we applied (\ref{glob}), and in (\ref{re})
we applied (\ref{l2emp}).
\hspace*{\fill} $\Box$

\subsection{Proof of Theorem \ref{ora-ent}}  \label{ora-ent-proof}

Let $f\in {\cal F}$ be the true density.
Let $\phi^0\in {\cal F}_{\delta}$.
Denote
$$
\zeta
=
C_1 \| \phi^0 - f\|_2^2 + C_2 n^{-1} \varrho^2(Q,{\cal F}_{\delta})
\log_e(\#{\cal F}_{\delta}) ,
$$
where
%\label{xi-choice}
$C_1$ is defined in (\ref{c1}) and $C_2$ is defined in (\ref{c2}).
%L_n = \log_e(\#{\cal F}_{\delta}) .
We have that
\begin{eqnarray}
\lefteqn{   \nonumber
E\| \hat{f} - f \|_2^2
}   \\  & = &   \nonumber
\int_0^{\infty} P \left( \| \hat{f} - f \|_2^2 > t \right) \, dt
\\ & \leq &   \nonumber
\zeta  +
\int_{\zeta}^{\infty}  P \left( \| \hat{f} - f \|_2^2 > t \right) \, dt
\\ & = &   \label{manyuse}
\zeta  +
C_2 n^{-1} \varrho^2(Q,{\cal F}_{\delta})
\int_{0}^{\infty}  P \left( \| \hat{f} - f \|_2^2
> C_2 n^{-1} \varrho^2(Q,{\cal F}_{\delta}) t + \zeta \right) \, dt .
\end{eqnarray}
Denote
$$
\tau_n
= C_{\tau} n^{-1} \varrho^2(Q,{\cal F}_{\delta})
\left( \log_e(\#{\cal F}_{\delta}) + t \right) ,
$$
where $C_{\tau}$ is defined in (\ref{ctau-2}).
Then,
\begin{eqnarray}
\lefteqn{   \nonumber
P \left( \| \hat{f} - f \|_2^2 > C_2 n^{-1}
\varrho^2(Q,{\cal F}_{\delta}) t + \zeta \right)
} \\ & = &  \nonumber
P \left( \| \hat{f} - f \|_2^2 >
C_1\|\phi^0 - f\|_2^2 + C_2 C_{\tau}^{-1}\tau_n \right)
\\ & = &   \nonumber
P \left( (1-2\xi)^{-1} \| \hat{f} - f \|_2^2 
\right. \\ && >   \left.  \nonumber
2\xi(1-2\xi)^{-1} \| \hat{f} - f \|_2^2 + C_1 \|\phi^0 - f\|_2^2
+ C_2 C_{\tau}^{-1}\tau_n
\right)
\\ & = &   \label{golu}
P \left( \| \hat{f} - f \|_2^2 >
2\xi\| \hat{f} - f \|_2^2 + (1+2\xi) \|\phi^0 - f\|_2^2 + \xi \tau_n
\right)   .
\end{eqnarray}
We have by Lemma \ref{oralemma},
$$
\left\| \hat{f} - f\right\|_2^2
\leq
\left\| \phi^0 - f\right\|_2^2
+
2 \nu_n[Q(\hat{f}-\phi^0)] .
$$
Denote
$$
w(\phi)
=
\| \phi - f \|_2^2 +  \|\phi^0 - f\|_2^2  + \tau_n/2 .
$$
Then we may continue (\ref{golu}) with
\begin{eqnarray}
\lefteqn{       \nonumber
P \left( \| \hat{f} - f \|_2^2 > C_2 n^{-1} \varrho^2(Q,{\cal F}_{\delta}) t
+ \zeta \right)
} \\ & & =     \nonumber
P \left( \nu_n[Q(\hat{f}-\phi^0)] >
\xi \| \hat{f} - f \|_2^2 +  \xi \|\phi^0 - f\|_2^2
+ \xi \tau_n /2  \right)
\\ & & =     \nonumber
P \left( \nu_n[Q(\hat{f}-\phi^0)] >  w(\hat{f}) \xi  \right)
\\ & & \leq      \nonumber
P \left( \max_{\phi\in{\cal F}_{\delta},\phi\neq \phi^0}
\frac{\nu_n[Q(\phi-\phi^0)]}{w(\phi)} >  \xi  \right) 
\\  & & \stackrel{def}{=}   \label{valistop}
P_{max}  .
\end{eqnarray}
We prove that
\begin{equation}  \label{expine-2}
P_{max}
\leq
\exp(-t) ,
\end{equation}
and this  proves the theorem, when we combine
(\ref{manyuse}) %, (\ref{golu}),
and (\ref{valistop}).

\paragraph{Proof of (\ref{expine-2})}
Denote
$$
{\cal G}
=
\left\{ \frac{Q(\phi - \phi^0)}{w(\phi)} \,
: \phi \in {\cal F}_{\delta}, \,\, \phi \neq \phi^0 \right\} .
$$
We have that
\begin{equation}  \label{sumout}
P_{max}
\leq
\sum_{g\in {\cal G}} P \left( \nu_n(g) > \xi \right) .
\end{equation}
Also,
$$ %\begin{equation}  \label{wala1}
w(\phi)
\geq
\frac{1}{2} \,  \left(
\left\| \phi - \phi^0 \right\|_2^2
+
\tau_n  \right)
\geq
\left\| \phi - \phi^0 \right\|_2  \tau_n^{1/2}
$$ %\end{equation}
and thus
\begin{equation}  \label{v-bound-white}
v_0
\stackrel{def}{=}
\max_{g\in {\cal G}} \|g\|_2^2
\leq
\frac{1}{\tau_n} \,
\max_{\phi\in {\cal F}_{\delta},\phi\neq \phi^0}
\frac{\|Q(\phi-\phi_0)\|_2^2}{\|\phi-\phi_0\|_2^2}
=
\frac{\varrho^2(Q,{\cal F}_{\delta})}{\tau_n} \,   .
\end{equation}

\bigskip  \noindent {\em Gaussian white noise.}
When $W\sim N(0,\sigma^2)$, then we have
$
P(W>\xi)  \leq 2^{-1} \exp\{-\xi^2/(2\sigma^2)\}
$
for $\xi>0$,
see for example \cite{Dudle99}, Proposition 2.2.1.
We have that $\nu_n(g) \sim N(0, n^{-1} \|g\|_2^2)$.
Thus,
$$
P \left( \nu_n(g) > \xi \right)
\leq
2^{-1} \exp\left\{-\,\frac{n\xi^2}{2v_0} \right\}
\leq
2^{-1} \exp\left\{
-\, \frac{n\tau_n\xi^2}{2\varrho^2(Q,{\cal F}_{\delta})}  \right\}  .
$$
Thus, denoting $C_{\xi}\stackrel{def}{=}\xi^2C_{\tau}/2$,
\begin{eqnarray*}
P_{max}
& \leq &
\#{\cal F}_{\delta} \cdot \exp\left\{
-\,\frac{n\tau_n\xi^2}{2\varrho^2(Q,{\cal F}_{\delta})}  \right\}
=
\#{\cal F}_{\delta} \cdot \exp\left\{ -C_{\xi}
[ \log_e(\#{\cal F}_{\delta}) + t ] \right\}
\\ & \leq &
\exp(-t) ,
\end{eqnarray*}
since $C_{\xi} \geq 1$ by the choice of $\xi$.

\bigskip \noindent {\em Density estimation.}
Denote
$v= \sup_{g\in {\cal G}} \mbox{Var}_f(g(Y_1))$,
and
$b= \sup_{g\in {\cal G}} \|g\|_{\infty}$.
We have that
\begin{equation}  \label{v-bound}
v
\leq
\|Af\|_{\infty} v_0
\leq
B_{\infty}\,\frac{\varrho^2(Q,{\cal F}_{\delta})}{\tau_n} \,,
\end{equation}
by (\ref{v-bound-white}).
Also,
$$ %\begin{equation}  \label{wala2}
w(\phi)
\geq
\frac{\tau_n}{2}
$$ %\end{equation}
and thus, because $\varrho(Q,{\cal F}_{\delta}) \geq 1$,
\begin{equation}  \label{b-bound}
b
\leq
2B_{\infty}' \, \frac{2}{\tau_n}\,
\leq
4 B_{\infty}' \, \frac{\varrho^2(Q,{\cal F}_{\delta})}{\tau_n}\, .
\end{equation}
Applying Bernstein's inequality,
applying (\ref{v-bound}) and (\ref{b-bound}),
\begin{eqnarray*}
P \left( \nu_n(g) > \xi \right)
& \leq &
\exp\left\{ \frac{-n\xi^2}{2(v+\xi b/3)} \right\}
\\ & \leq &
\exp\left\{ \frac{-n\xi^2\tau_n}{2\varrho^2(Q,{\cal F}_{\delta})
(B_{\infty}+4B_{\infty}'\xi/3)} \right\} .
\end{eqnarray*}
Continuing from (\ref{sumout}),
\begin{eqnarray*}
P_{max}
&\leq&
\#{\cal F}_{\delta} \cdot \exp\left\{ \frac{-n\xi^2\tau_n}{2\varrho^2(Q,{\cal F}_{\delta})
(B_{\infty}+4B_{\infty}'\xi/3)} \right\}
\\ &=&
\#{\cal F}_{\delta} \cdot \exp\left\{ -C_{\xi}
[ \log_e(\#{\cal F}_{\delta}) + t ] \right\}
\\ &\leq &
\exp(-t)
\end{eqnarray*}
where
$$
C_{\xi} \stackrel{def}{=}
\frac{\xi^2 C_{\tau}}{2(B_{\infty}+4B_{\infty}'\xi/3)} \,,
$$
and $C_{\xi} \geq 1$ by the choice of $\xi$.
%$$
%C_{\tau} \xi^2 - 8B_{\infty}\xi/3 - 2B_{\infty} \geq 0 .
%$$
We have proved (\ref{expine-2}) and thus the theorem.
\hspace*{\fill}   $\Box$  %$\square$ %$\Box$

%\subsection{Proof of Theorem~\ref{packing-low}}
\subsection{Proof of Theorem~\ref{lowtheo}}    \label{lowtheo-proof}

%We prove a more general proposition than Theorem~\ref{packing-low}.
%Theorem~\ref{packing-low} follows then as a corollary.

To prove Theorem~\ref{lowtheo}
we follow the approach of \cite{HasmiIbrag90}.
We start with a useful lemma.

\begin{lemma}  \label{summari}
Let ${\cal D} \subset {\cal F}$ be a finite set for which
\begin{equation} \label{separe}
\min \{  \|f-g\|_2  :
f,g \in {\cal D}  ,
\,\,\,
f \neq g
\}   \geq \delta
\end{equation}
where $\delta >0$.
Assume that for some $f_0 \in {\cal D}$,
and for all $f \in {\cal D} \setminus \{ f_0 \}$,
\begin{equation}  \label{likfraprod}
P_{Af}^{(n)} \left( \frac{dP_{Af_0}^{(n)}}{dP_{Af}^{(n)}} \leq \tau \right)
\leq \alpha  ,
\end{equation}
where
$0<\alpha<1$, $\tau>0$, and
in the density estimation model
$P_{Af}^{(n)}$ is the product measure corresponding to density $Af$,
and in the Gaussian white noise model
$P_{Af}^{(n)}$ is the measure of process $Y_n$ in (\ref{yn}).
Then,
$$
\inf_{\hat{f}} \sup_{f\in {\cal F}} E_{Af} \|f-\hat{f}\|_2^2
\geq
\frac{\delta^2}{4} \,
(1-\alpha)\, \frac{\tau (N_{\delta}-1)}{1+\tau (N_{\delta}-1)} \, ,
$$
where
$N_{\delta} = \#{\cal D} \geq 2$
and the infimum is taken over all estimators
(either in the density estimation model or in the Gaussian white noise model).
\end{lemma}

{\em Proof.}
Let $f_n: {\bf R}^d \to {\bf R}$ be an estimator of $f$.
Define a random variable $\hat{\theta}$ taking values in ${\cal D}$:
$$
\hat{\theta} = \mbox{argmin}_{f \in {\cal D}}  \|f_n-f\|_2  .
$$
Note that by (\ref{separe}),
$$
\hat{\theta} \neq f \in {\cal D}  \Rightarrow \|f_n-f\|_2 \geq \delta/2 ,
$$
since $\hat{\theta} \neq f$ for an $f \in {\cal D}$ implies that
$f_n$ is closer to some other $g \in {\cal D}$ than to $f$.
Then, applying also Markov's inequality,
\begin{eqnarray*}
\sup_{f \in {\cal F}}  E_{Af} \|f_n-f\|_2^2
& \geq &
\max_{f \in {\cal D}}  E_{Af} \|f_n-f\|_2^2
\\
& \geq &
\frac{\delta^2}{4} \, \max_{f \in {\cal D}} P_{Af}^{(n)}
\left( \|f_n-f\|_2^2  \geq \delta^2/4 \right)
\\
& \geq &
\frac{\delta^2}{4} \, \max_{f \in {\cal D}}
P_{Af}^{(n)} ( \hat{\theta} \neq f ) .
\end{eqnarray*}
The lemma follows by an application of Lemma~\ref{tsybalemma} below.
%Note that we have used the typical notation $P_f ( \hat{\theta} \neq f )$
%instead of $P_f^{(n)} ( \hat{\theta} \neq f )$.
\hspace*{\fill} $\Box$

\begin{lemma}  \label{tsybalemma}
(\cite{Tsyba98}, Theorem 6.)
Let $\hat{\theta}$ be a random variable taking values on
a finite set ${\mathbb P}$ of probability measures.
Denote $\#{\mathbb P} = N$ and assume $N\geq 2$.
Let $\tau>0$ and $0<\alpha<1$. Let for some $P_0 \in {\mathbb P}$
and for all $P \in {\mathbb P} \setminus \{ P_0 \}$,
\begin{equation}  \label{likfra}
P \left( \frac{dP_0}{dP} \leq \tau \right) \leq \alpha  .
\end{equation}
Then
$$
\max_{P \in {\mathbb P}} P ( \hat{\theta} \neq P )
\geq
(1-\alpha)\, \frac{\tau (N-1)}{1+\tau (N-1)} \, .
$$
\end{lemma}

\paragraph{Proof of Theorem~\ref{lowtheo}}

For $f,f_0\in {\cal D}_{\psi_n}$, $f\neq f_0$,
\begin{eqnarray}  \nonumber
\lefteqn{
P_{Af}^{(n)} \left( \frac{dP_{Af_0}^{(n)}}{dP_{Af}^{(n)}} \leq \tau \right)
} \\ & \leq &     \label{eaqi}
\left( \log \tau^{-1} \right)^{-1} D_K^2(P_{Af}^{(n)},P_{Af_0}^{(n)})
\\ & = &    \label{taqi}
\left\{  \begin{array}{ll}
\left( \log \tau^{-1} \right)^{-1} n D_K^2(Af,Af_0),
& \mbox{ density estimation }
\\
\left( \log \tau^{-1} \right)^{-1} \frac{n}{2} \,  \|Af-Af_0\|_2^2,
& \mbox{Gaussian white noise, }
\end{array} \right.
\end{eqnarray}
where
in (\ref{eaqi}) we applied Markov's inequality,
and
in (\ref{taqi}) we applied
for the Gaussian white noise model the fact that under $P_{Af}^{(n)}$,
$$
\frac{dP_{Af}^{(n)}}{dP_{Af_0}^{(n)}}
=
\exp\left\{ n^{1/2} \sigma Z + n\sigma^2/2 \right\} ,
$$
where
$Z\sim N(0,1)$ and $\sigma= \|Af-Af_0\|_2$.
When we choose
$$
\tau = \tau_n = \exp \left\{ -\alpha^{-1} n [C_1 \varrho_K(A,{\cal D}_{\psi_n} ) \psi_n]^2
\right\},
$$
for $0< \alpha <1$, then
applying assumption (\ref{klassu}),
\begin{eqnarray}
P_{Af}^{(n)} \left( \frac{dP_{Af_0}^{(n)}}{dP_{Af}^{(n)}} \leq \tau \right)
& \leq &  \nonumber
\left( \log \tau^{-1} \right)^{-1} n  \varrho_K^2(A,{\cal D}_{\psi_n} )  \|f-f_0\|_2^2
\\ & \leq &  \nonumber
\left( \log \tau^{-1} \right)^{-1} n [ \varrho_K(A,{\cal D}_{\psi_n} ) C_1\psi_n ]^2
\\ & = &   \label{naqi}
\alpha .
\end{eqnarray}
Applying Lemma \ref{summari}, assumption (\ref{l2assu}), and (\ref{naqi})
we get the lower bound
\begin{equation}  \label{fiinis}
\inf_{\hat{f}} \sup_{f\in {\cal D}_{\psi_n}} \|f-\hat{f}\|_2^2
\geq
\frac{(C_0\psi_n)^2}{4} \,
(1-\alpha) \,
\frac{\tau_n (N_{\psi_n}-1)}{1+\tau_n (N_{\psi_n}-1)} \, ,
\end{equation}
where $N_{\psi_n} = \# {\cal D}_{\psi_n}$.
Let $n$ be so large that
$\log_e N_{\psi_n} \geq C_2^2 n \varrho_K^2(A,{\cal D}_{\psi_n} )  \psi_n^2$,
where $C_2>C_1$.
This is possible by (\ref{minimaxeq}).
Then,
\begin{eqnarray*}
\tau_nN_{\psi_n}
& = &
\exp \left\{ \log_e N_{\psi_n} -\alpha^{-1} n [C_1 \varrho_K(A,{\cal D}_{\psi_n} )
\psi_n]^2 \right\}
\\ & \geq &
\exp \left\{ n \varrho_K^2(A,{\cal D}_{\psi_n} ) \psi_n^2 [C_2^2-\alpha^{-1}C_1^2] \right\}
\to \infty
\end{eqnarray*}
as $n\to \infty$,
where
we apply (\ref{epsi-infi}) and
choose $\alpha$ so that
$C_2^2-\alpha^{-1}C_1^2 >0$,
that is,
$(C_1/C_2)^2 < \alpha < 1$.
Then
$$
\lim_{n\to \infty}
\frac{\tau_n (N_{\psi_n}-1)}{1+\tau_n (N_{\psi_n}-1)}
=1
$$
and the theorem follows from (\ref{fiinis}).
\hspace*{\fill} $\Box$

\subsection{Proof of Theorem~\ref{chain-white}}  \label{proof-chain-whitenoise}

Denote
$$
\zeta = C_1 \epsilon + C_2 \psi_n^2  ,
$$
where $C_1=(1-2\xi)^{-1}$, $C_2=1-2\xi$,
$0< \xi \leq (3-\sqrt{5})/4$.
We have that
\begin{eqnarray}
E\| \hat{f} - f \|_2^2
& = &       \nonumber
\int_0^{\infty} P \left( \| \hat{f} - f \|_2^2 > t \right) \, dt
\\ & \leq &    \nonumber
\zeta +
\int_{\zeta}^{\infty}
P \left( \| \hat{f} - f \|_2^2 > t \right) \, dt
\\ & = &    \label{zappi}
\zeta +
C_2 \psi_n^2
\int_{0}^{\infty}  P \left( \| \hat{f} - f \|_2^2
> C_2 \psi_n^2 t + \zeta  \right) \, dt .
\end{eqnarray}
Denote
$$
\tau_n = C_{\tau} \psi_n^2( 1+t) ,
\qquad
C_{\tau} = \xi^{-1} (1-2\xi)^2 .
$$
Then,
\begin{eqnarray}
\lefteqn{   \nonumber
P \left( \| \hat{f} - f \|_2^2 > C_2 \psi_n^2 t + \zeta \right)
} \\ & = &  \nonumber
P \left( \| \hat{f} - f \|_2^2 > C_2C_{\tau}^{-1}\tau_n + C_1\epsilon \right)
\\ & = &   \nonumber
P \left( (1-2\xi)^{-1} \| \hat{f} - f \|_2^2 
\right. \\ && >  \nonumber  \left.
2\xi(1-2\xi)^{-1} \| \hat{f} - f \|_2^2 +   C_2C_{\tau}^{-1} \tau_n
+ C_1 \epsilon
\right)
\\ & = &   \label{golu-2}
P \left( \| \hat{f} - f \|_2^2 >
2\xi\| \hat{f} - f \|_2^2 + \xi \tau_n + \epsilon
\right)   .
\end{eqnarray}
We have by Lemma \ref{oralemma}, choosing $f^0=f$,
$$
\left\| \hat{f} - f\right\|_2^2
\leq
2 \nu_n[Q(\hat{f}-f)] + \epsilon .
$$
Denote
$$
w(g)
=
\| g - f \|_2^2   + \tau_n/2 .
$$
Then we may continue (\ref{golu-2}) with
\begin{eqnarray}
\lefteqn{      \nonumber
P \left( \| \hat{f} - f \|_2^2 > C_2 \psi_n^2 t + \zeta \right)
} \\ & & \leq    \nonumber
P \left( \nu_n[Q(\hat{f}-f)] >
\xi \| \hat{f} - f \|_2^2  + \xi \tau_n /2  \right)
\\ & & =     \nonumber
P \left( \nu_n[Q(\hat{f}-f)] >  w(\hat{f}) \xi  \right)
\\ & & \leq     \nonumber
P \left( \sup_{g\in{\cal F}}
\frac{\nu_n[Q(g-f)]}{w(g)} >  \xi  \right) 
\\  & & \stackrel{def}{=}   \label{zuppi}
P_{sup}  .
\end{eqnarray}
We prove that
\begin{equation}  \label{expine-2-2}
P_{sup}
\leq
\exp(-t \cdot \log_e 2) ,
\end{equation}
and this proves the theorem, when we combine (\ref{zappi}) and (\ref{zuppi}).

\paragraph{Proof of (\ref{expine-2-2})}
We use the peeling device, see for example
\cite{vanGe00}, page 69.
Denote
$$
a_0 =  \tau_n/2 ,
\qquad
a_j= 2^{2j} a_0 ,
\qquad
b_j = 2^2a_j,
\qquad
j=0,1,\ldots
$$
Let ${\cal G}_j$ be the set of functions
$$
{\cal G}_j
=
\left\{ g \in {\cal F} :
a_j \leq w(g) < b_j \right\} ,
\qquad
j=0,1,\ldots
$$
and
$$
{\cal F}_j
=
\left\{ g \in {\cal F} :
\|g-f\|_2^2 < b_j \right\} ,
\qquad
j=0,1,\ldots
$$
We have that
$$
{\cal F}
=
\{ g \in {\cal F} : w(g) \geq a_0 \}
=
\bigcup_{j=0}^{\infty} {\cal G}_j .
$$
Thus,
\begin{eqnarray}
P_{sup}
& \leq &   \nonumber
\sum_{j=0}^{\infty}
P\left( \sup_{g\in {\cal G}_j} \frac{\nu_n[Q(g-f)]}{w(g)} > \xi \right)
\\ & \leq &    \nonumber
\sum_{j=0}^{\infty}
P\left( \sup_{g\in {\cal F},w(g)<b_j} \nu_n[Q(g-f)] > \xi a_j \right)
\\ & \leq &  \label{sandw-2-2}
\sum_{j=0}^{\infty}
P\left( \sup_{g\in {\cal F}_j} \nu_n[Q(g-f)] > \xi a_j \right)  .
\end{eqnarray}
By Assumption~4 of Theorem~\ref{chain-white},
$\tilde{G}(\psi_n)=24\sqrt{2} G(\psi_n)$,
where $\tilde{G}$ is defined in (\ref{xicond-white}),
for sufficiently large $n$.
Thus, by the choice of $C=\xi^{-1}4\cdot 24\sqrt{2}$ in (\ref{rateeq}),
$$
\psi_n^2  \geq n^{-1/2} \xi^{-1}4 \tilde{G}(\psi_n)  .
$$
By the choice of $\xi$ we have that $C_{\tau} \geq 2$, and thus
$a_0 = C_{\tau}\psi_n^2 (1+t)/2 \geq  \psi_n^2$.
Since $G(\delta)/\delta^2$ is decreasing,
by Assumption~2 of Theorem~\ref{chain-white},
then $\tilde{G}(\delta)/\delta^2$ is decreasing,
and
$$
\xi n^{1/2} /4
\geq \tilde{G}(\psi_n)/\psi_n^2
\geq \tilde{G}(a_0^{1/2})/a_0
\geq \tilde{G}(b_j^{1/2})/b_j ,
$$
that is,
\begin{equation}  \label{cjappu-2-2}
\xi a_j = \xi b_j/4  \geq n^{-1/2} \tilde{G}(b_j^{1/2})  .
\end{equation}
We may apply Lemma~\ref{telehoef-white} 
given in Appendix \ref{appe-white},
with (\ref{cjappu-2-2}) to get
\begin{eqnarray}
\lefteqn{      \nonumber
P\left( \sup_{g\in {\cal F}_j} \nu_n[Q(g-f)] > \xi a_j \right)
} \\ &\leq &  \label{diffestep}
\exp\left\{ -\,\frac{ n(\xi a_j)^2C'}{c^2 b_j^{1-a}} \right\}
\\ & \leq &   \nonumber
\exp\left\{ - C'' 2^{2j(a+1)} n\psi_n^{2(1+a)} (1+t)^{1+a} \right\}
\\ & \leq &    \label{repo}
\exp\left\{ - C'' (j+1) n\psi_n^{2(1+a)} (1+t)^{1+a} \right\}  ,
\end{eqnarray}
where
$C'' = C' c^{-2} \xi^2 2^{2(a-1)}(C_{\tau}/2)^{1+a}$,
and we used the facts
$a_j^2/b_j^{1-a}
=2^{2(a-1)}a_j^{1+a}
=2^{2(a-1)} (2^{2j}a_0)^{1+a}
=2^{2(a-1)} [2^{2j}C_{\tau}\psi_n^2(1+t)/2]^{1+a}$
and
$2^{2j(a+1)} \geq j+1$.
When $0 \leq b \leq 1/2$, then
$\sum_{j=0}^{\infty} b^{j+1} = \sum_{j=1}^{\infty} b^j = b/(1-b) \leq 2b$.
When
$n\psi_n^{2(1+a)} \geq (\log_e2) / C''$,
then
$\exp\{ - C'' n\psi_n^{2(1+a)} (1+t)^{1+a} \} \leq 1/2$,
and
we combine (\ref{sandw-2-2}) and (\ref{repo}) to get the upper bound
\begin{eqnarray*}
2 \exp\left\{ - C'' n\psi_n^{2(1+a)} (1+t)^{1+a} \right\}
& \leq &
2 \exp\left\{ - C'' n\psi_n^{2(1+a)} (1+t) \right\}
\\ & \leq &
\exp\left\{ - t \log_e2  \right\}  .
\end{eqnarray*}
We have proved (\ref{expine-2-2}) and thus
we have proved Theorem~\ref{chain-white}
up to proving Lemma~\ref{telehoef-white}, 
which is done in Appendix \ref{appe-white}.

\subsection{Proof of Theorem~\ref{chain-density}}  \label{proof-chain-density}

The proof goes similarly as the Proof of Theorem~\ref{chain-white}
until step (\ref{diffestep}).
At this step we apply Lemma~\ref{telehoef-density},
given in Appendix \ref{appe-empiri},
to get
\begin{eqnarray*}
\lefteqn{
P\left( \sup_{g\in {\cal F}_j} \nu_n[Q(g-f)] > \xi a_j \right)
} \\ & \leq &
\exp\left\{ -\,\frac{ n(\xi a_j)^2C'}{c^2 b_j^{1-a}} \right\}
+
2 \#{\cal G}_{B_2} \exp\left\{ -\, \frac{1}{12} \,
\frac{ n (\xi a_j)^2}{B_{\infty}c^2b_j^{1-a}
+ 2\xi a_j B_{\infty}'/9} \right\}  .
\end{eqnarray*}
The first term in the right hand side is handled similarly as
in the Proof of Theorem~\ref{chain-white}.
For the second term in the right hand side we have,
for sufficiently large $n$,
\begin{eqnarray*}
\exp\left\{ -\, \frac{1}{12} \,
\frac{ n (\xi a_j)^2}{B_{\infty}c^2b_j^{1-a}
+ 2\xi a_j B_{\infty}'/9} \right\}
& = &
\exp\left\{ -\, \frac{1}{12} \,
\frac{ n \xi^2 a_j}{B_{\infty}c^2 a_0^{-a}
+ 2\xi B_{\infty}'/9} \right\}
\\ & \leq &
\exp\left\{ -\, \frac{1}{12} \,
\frac{ n \xi^2 a_j a_0^a }{B_{\infty}c^2
+ 2\xi B_{\infty}'/9} \right\}
\\ & = &
\exp\left\{ - n \psi_n^{2(1+a)} 2^{2j} (1+t)^{1+a} C''  \right\},
\end{eqnarray*}
since $a_j^{-a} = (2^{2j}a_0)^{-a} \leq a_0^{-a}$
and $a_0^{-a} \geq 1$ for sufficiently large $n$,
and we denote
$C' = \xi^2 C_{\tau}^{1+a}/[2^{1+a}12(B_{\infty}c^2 + 2\xi B_{\infty}'/9)]$.
The proof is finished similarly as the proof of Theorem~\ref{chain-white}.

\subsection{Proof of Theorem \ref{ora-add}}  \label{ora-add-proof}

We proceed similarly as in the proof of Theorem~\ref{ora-ent}. 
Choose $f_{\delta} \in {\cal F}_\delta$ such that 
$\|f - f_{\delta}\|_2 \leq \delta$, where $f$ is the underlying function in 
${\cal F}$.
Choose $\xi< 1/2$ and put $\zeta = \zeta_1 + \zeta_2$ with 
$\zeta_1 = (1-2 \xi)^{-1} (1+2 \xi) \|f - f_{\delta}\|_2 ^2$, 
$\zeta_2 = \kappa n ^{-1} \sum_{j=1}^p \rho_j^2 \lambda_j$ and 
$\kappa = 4 c^{-1} \xi^{-1} (1-2 \xi)^{-1}$.
We have that
\begin{eqnarray}
E \left( \|\hat f - f\|_2^2 \right) 
&\leq& 
\zeta 
+ \int_ \zeta^{\infty} P\left( \|\hat f - f\|_2^2 > t \right)\, dt 
\nonumber
\\
\label{geo2} 
&\leq& 
\zeta 
+ 
\int_ 0^{\infty} P\left( \|\hat f - f\|_2^2 > t+ \zeta \right)\, dt.
\end{eqnarray}
For the integrand of the second term we have
that
\begin{eqnarray*}
\lefteqn{
P\left(\|\hat f - f\|_2^2 > t+ \zeta\right) 
} \\ & = & 
P\left( (1-2 \xi)^{-1}\|\hat f - f\|_2^2
 > 2 \xi (1-2 \xi)^{-1}\|\hat f - f\|_2^2+ t+ \zeta \right) 
\\ & = & 
P \left(\|\hat f - f\|_2^2 > 
2 \xi \|\hat f - f\|_2^2+ (1-2 \xi)t+ (1-2 \xi)\zeta \right)  .
\end{eqnarray*}
We now use Lemma~\ref{oralemma}. This gives
\begin{eqnarray*}
\|\hat f - f\|_2^2 
\leq 
\| f - f_{\delta}\|_2^2+ 2 \nu_n\left(Q(\hat f - f_{\delta}) \right).
\end{eqnarray*}
Together with the last equalities this gives
\begin{eqnarray*}
\lefteqn{
P\left(\|\hat f - f\|_2^2 > t+ \zeta \right) 
} \\ & \leq &
P\left( \| f - f_{\delta}\|_2^2+ 2 \nu_n \left(Q(\hat f - f_{\delta}) \right) 
%\right. \\ && \qquad \qquad \qquad \left.
> 2 \xi \|\hat f - f\|_2^2+ (1-2 \xi)(t+\zeta)  \right)
\\ & = &
P\left( \nu_n\left(Q(\hat f - f_{\delta})\right) 
\right. \\ && \qquad \qquad \qquad \left.
>  \xi \|\hat f - f\|_2^2+ \xi \| f - f_{\delta}\|_2^2 
+ 2^{-1}(1-2 \xi)(t+\zeta_2)  \right)
\\ & \leq & 
P\left( \nu_n\left(Q(\hat f - f_{\delta})\right) 
>  2^{-1} \xi \|\hat f - f_{\delta}\|_2^2 
+ 2^{-1}(1-2 \xi)(t + \zeta_2)  \right).
\end{eqnarray*}
Now, put 
$w_j=\rho_j / \sum_{l=1}^p \rho_l$ 
and decompose 
$f_{\delta} = f_{\delta,1}+\cdots+f_{\delta,p}$ and 
$\hat f= \hat f_1 +\cdots+\hat f_p$ with 
$f_{\delta,j}, \hat f_j \in {\cal F}_{j,\delta}$. 
Using assumption (\ref{geo1}) we get with 
$\beta_j = 2^{-1}(1-2 \xi) (w_jt + \kappa n^{-1} \rho_j^2 \lambda_j)$,
\begin{eqnarray*}
\lefteqn{
P\left(\|\hat f - f\|_2^2 > t+ \zeta\right) 
} \\ & \leq & 
P\left( \sum_{j=1}^p \nu_n\left(Q(\hat f_j - f_{\delta,j})\right) 
> 2 ^{-1} \xi c \sum_{j=1}^p\|\hat f_j  - f_{\delta,j}\|_2^2 
+ \sum_{j=1}^p \beta_j \right)
\\ & \leq &
\sum_{j=1}^p  P\left( \nu_n\left(Q(\hat f_j - f_{\delta,j}) \right) 
> 2 ^{-1} \xi c \|\hat f_j  - f_{\delta,j}\|_2^2 + \beta_j\right)
\\ & \leq &
\sum_{j=1}^p \sum_{g_j \in {\cal F}_{j,\delta}} 
P\left(   \nu_n\left(Q(g_j - f_{\delta,j}) \right) 
> 2 ^{-1} \xi c \|g_j  - f_{\delta,j}\|_2^2  + \beta_j \right).
\end{eqnarray*}
We now use 
$$ 
P\left(\nu_n(h) > \xi \right) 
\leq 
2^{-1} \exp \left (-{n \xi^2 \over 2 \|h\|^2_2} \right ),
$$
compare to the proof of Theorem~\ref{ora-ent}. This gives
\begin{eqnarray*}
\lefteqn{
P\left(\|\hat f - f\|_2^2 > t+ \zeta \right) 
} \\ & \leq &
\sum_{j=1}^p \sum_{g_j \in {\cal F}_{j,\delta}} 2 ^{-1} 
\exp \left [ - {n (2 ^{-1} \xi c \|g_j  - f_{\delta,j}\|_2^2 + \beta_j)^2 
\over 2 \|Q(g_j - f_{\delta,j})\|^2} \right ]
\\ & \leq &
\sum_{j=1}^p \sum_{g_j \in {\cal F}_{j,\delta}} 2 ^{-1} 
\exp \left [ - { n \xi c \|g_j  - f_{\delta,j}\|_2^2  \beta_j 
\over 2 \|Q(g_j - f_{\delta,j})\|^2} \right ]
%\\ && \qquad  \leq
%\sum_{j=1}^p \sum_{g_j \in {\cal F}_{j,\delta}} 2 ^{-1} 
%\exp \left [ - { n \xi c  \|g_j  - f_{\delta,j}\|_2^2 \beta_j \over 2 \rho_j^2 
%\|g_j - f_{\delta,j}\|_2^2} \right ] 
\\ & \leq &
\sum_{j=1}^p \exp (\lambda_j) 2 ^{-1} 
\exp \left [ - {  n\xi c  \beta_j \over 2 \rho_j^2 } \right ]
\\ & = & 
\sum_{j=1}^p  2 ^{-1} 
\exp \left [ - n \xi c 4^{-1} (1-2 \xi) w_j \rho_j^{-2}t \right ].
\end{eqnarray*}
By plugging this into (\ref{geo2}) we get
\begin{eqnarray*}
E \left( \|\hat f - f\|_2^2\right) 
&\leq& 
 \zeta + \sum_{j=1}^p \int_ 0^{\infty} 
\exp \left [ - n \xi c 4^{-1} (1-2 \xi) w_j \rho_j^{-2}t \right ] \, dt \\
&\leq&
\zeta + \sum_{j=1}^p   n^{-1} 4 [\xi c  (1-2 \xi)w_j]^{-1}  \rho_j^{2}\\
&=&
\zeta +    n^{-1} 4 [\xi c  (1-2 \xi)]^{-1} 
\left(\sum_{j=1}^p \rho_j\right)^{2} .
\end{eqnarray*}
Choosing $\xi = 4 ^{-1}$ gives the statement of Theorem~\ref{thm4}.

\section*{Acknowledgment}

We would like to thank referees for suggesting improvements and
pointing out errors.
Writing of this article was financed by
Deutsche Forschungsgemeinschaft under project MA1026/8-1.

\bibliographystyle{agsm}
\bibliography{../viite2}

\appendix

\section{Ellipsoids}  \label{ellipsoid-appendix}

The ellipsoid has been defined in (\ref{sobolev-ellipsoid})
and we assume that the $a_j$ satisfy (\ref{abound}).
We make the calculations now in the one dimensional case.

\subsection{$\delta$-net}

We shall construct a $\delta$-net ${\Theta}_{\delta}$
for the ellipsoid in (\ref{sobolev-ellipsoid}).
The construction is similar to the construction of
\cite{KolmoTikho61}.
Let
\begin{equation}  \label{M-definition-II}
M = [ ( C_1^{-1} 2^{1/2} L\delta^{-1})^{1/s}].
\end{equation}
Let $\Theta_{\delta}(M)$ be a $\delta/2$-net of
$$
E_M
=
\left\{ (\theta_j)_{j\in \{ 1,\ldots ,M\}}  :
\sum_{j=1}^{M}  a_j^2  \theta_j^2
\leq L^2 \right\}  .
$$
We can choose $\Theta_{\delta}(M)$ in such a way that its cardinality
satisfies
$$
\# \Theta_{\delta}(M)
\leq
C \, \frac{\mbox{volume}(E_M)}{\mbox{volume}(B_{\delta}^{(M)})} \, ,
$$
where
$B_{\delta}^{(M)}$
is a ball of radius $\delta$ in the $M$-dimensional Euclidean space.
Define the $\delta$-net by
\begin{eqnarray} \nonumber
\Theta_{\delta}
=
\left\{  (\theta_j)_{j\in \{ 1,\ldots ,\infty\}}
:
 (\theta_j)_{j\in \{ 1,\ldots ,M\}}   \in \Theta_{\delta}(M),
\, \theta_j=0, \, \mbox{ for } j \geq M+1  \right\}  .
\label{delta-net-definition-II}
\end{eqnarray}

\bigskip \noindent
{\em ($\delta$-net property.)}
We proof that $\Theta_{\delta}$
is a $\delta$-net of the ellipsoid $\Theta$.
For each $\theta\in \Theta$ there is
$\theta_{\delta} \in \Theta_{\delta}$
such that $\|\theta - \theta_{\delta}\|_{l_2} \leq \delta$.
Indeed, let $\theta\in \Theta$.
Let $\theta_{\delta} \in \Theta_{\delta}$ be such that
$$
\sum_{j=1}^{M}
(\theta_j - \theta_{\delta,j} )^2  \leq \delta^2/2  .
$$
Then
$$
\| \theta - \theta_{\delta} \|_{l_2}^2
=
\sum_{j=1}^{M}
(\theta_j - \theta_{\delta,j} )^2
+
\sum_{j=M+1}^{\infty}  \theta_j^2
\leq
\delta^2
$$
where we used the fact
\begin{equation}  \label{fact}
\sum_{j=M+1}^{\infty}
\theta_j^2
\leq
C_1^{-2} \cdot M^{-2s}
\sum_{j=M+1}^{\infty} a_j^2 \theta_j^2
\leq
C_1^{-2} M^{-2s} L^2
\leq
\delta^2/2 ,
\end{equation}
because, when $j \notin \{ 1,\ldots ,M\}$, then
$$
a_j^{-2} \leq C_1^{-2} \cdot j^{-2s} \leq C_1^{-2} \cdot M^{-2s}
\leq
\delta^2/(2L^2)  .
$$

\bigskip \noindent
{\em (Cardinality.)}
We  prove that
$$ %\begin{equation}  \label{ellipsoid-cardinality-I}
\log( \# \Theta_{\delta})
\leq
C \delta^{-1/s} .
$$ %\end{equation}
We have that
$$
\mbox{volume}(E^{(M)})
=
C_M  \cdot  L^{M}
\prod_{j=1}^{M} a_j^{-1}
$$
and
$$
\mbox{volume}(B_{\delta}^{(M)})
=
C_M  \cdot \delta^{M} ,
$$
where $C_M$ is the volume of the unit ball in the $M$ dimensional
Euclidean space.
Thus the cardinality of $\Theta_{\delta}$ satisfies
$$
\# \Theta_{\delta}
=
\# \Theta_{\delta}(M)
\leq
C\, \frac{L^{M} \prod_{j=1}^{M} a_j^{-1}}{\delta^{M}} \, .
$$
We have that
$$
%C' (M!)^{-s} = C'  \prod_{k=1}^{M} k^{-s} \leq
\prod_{j=1}^{M}  a_j^{-1}
\leq
C \prod_{j=1}^{M} j^{-s}
=
C (M!)^{-s} .
$$
Applying \cite{Felle68}, pp. 50-53, we get
$$
M!
>
M^{M+1/2}e^{-M}  .
%< 2 M^{M+1/2}e^{-M}  .
$$
Thus
\begin{eqnarray}
\lefteqn{   \nonumber
\log ( \# \Theta_{\delta} )
} \\ & \leq &  \nonumber
M \log(L) -s\log(M!) + M \log(\delta^{-1})  + C
\\ & \leq &   \nonumber
M \log(L) - s(M+1/2) \log M + sM
+ M \log(\delta^{-1})  +C
\\ & \leq &   \nonumber
M (\log(L)+s) - sM \log M + M \log(\delta^{-1})  +C
\\ & \leq &   \nonumber
M (\log(L)+s + C') + C
\\ & \leq &   \label{cardi-final}
\delta^{-1/s} C'' + C  ,
\end{eqnarray}
since $M = C''' \delta^{-1/s}$.

\subsection{$\delta$-packing set}

For a fixed sequence $\theta^*$ with
$\sum_{j=0}^{\infty} a_j^2 {\theta_j^*}^2 = L^* < L$ let
$\Theta^*_{\delta}(M)$ be a $\delta$-packing set of
$$
E^*_M
=
\left\{ (\theta_j)_{j\in \{ M^*,\ldots ,M\}}  :
\sum_{j=M^*}^{M}  a_j^2  \theta_j^2 \leq (L-L^*)^2
\right\} .
$$
Here, $M^* = [ M/2 ]$.
We can choose $\Theta^*_{\delta}(M)$
in such a way that its cardinality satisfies
\begin{equation}  \label{ellipsoid-cardinality2}
\log( \# \Theta^*_{\delta}(M)) \geq C^* \delta^{-1/s} .
\end{equation}
Define
\begin{eqnarray}
\Theta^*_{\delta}
=     \nonumber
\left\{  (\theta_j)_{j\in \{ 0,\ldots ,\infty\}}
\right.  & : &
(\theta_j- \theta_j^*)_{j\in \{ M^*,\ldots ,M\}} \in \Theta^*_{\delta}(M),
\\ &&   \left.
\, \theta_j=\theta_j^*, \, \mbox{ for } j \notin
\{ M^*,\ldots , M\} \right\} .
\label{delta-net-definition2}
\end{eqnarray}
The bound (\ref{ellipsoid-cardinality2}) follows
similarly as the upper bound (\ref{cardi-final}). In the
white noise case one can use this construction with $\theta^*=0$ and $L^*
=0$. In the density case another choice of $\theta^*$ may be appropriate
to ensure that the functions in ${\cal D}_{\delta}$ are bounded
from above and from below. This would allow to use the bound
(\ref{klbounds}) to carry over bounds on Hilbert norms to
corresponding bounds on Kullback-Leibler distances.
Note also that a similar calculation as in (\ref{fact})
shows that
for $\theta,\theta' \in \Theta_{\delta}^*$,
\begin{equation}  \label{klassu-holds}
\| \theta - \theta' \|_{l_2}^2
=
\sum_{i=M^*}^{M}  ( \theta_i -\theta_i' )^2
=
\sum_{i=M^*}^{\infty}  ( \theta_i -\theta_i' )^2
\leq
C \delta^2  .
\end{equation}

\section{Lemmas related to empirical process theory}

\subsection{Gaussian white noise}  \label{appe-white}

Lemma \ref{telehoef-white} gives an exponential tail bound for the
Gaussian white noise model.

\begin{lemma}  \label{telehoef-white}
Let $\nu_n$ be the centered empirical operator
of a Gaussian white noise process.
Operator $\nu_n$ is defined in (\ref{nun-def}).
%Let $Y_n$ be the Gaussian white noise process with signal function $Af$,
Let ${\cal G}\subset L_2({\bf R}^d)$
be such that $\sup_{g\in{\cal G}} \|g\|_2 \leq R$
and denote with ${\cal G}_{\delta}$ a $\delta$-net of ${\cal G}$,
$\delta>0$.
Assume that $\delta \mapsto \varrho(Q,{\cal G}_{\delta}) \sqrt{\log_e(\#{\cal G}_{\delta})}$
is decreasing on $(0,R]$,
where $\varrho(Q,{\cal G}_{\delta})$ is defined in (\ref{t-dense})
and
assume that the entropy integral $G(R)$
defined in (\ref{ent-int}) is finite.
Assume that $\varrho(Q,{\cal G}_{\delta}) = c\delta^{-a}$,
where $0\leq a <1$ and $c>0$.
Then for all
\begin{equation}  \label{xicond-white}
\xi \geq  n^{-1/2}\, \tilde{G}(R),
\qquad
\tilde{G}(R)
=
\max\left\{ 24\sqrt{2} G(R)  ,
c R^{1-a} \sqrt{\log_e 2 / C'} \right\}
\end{equation}
where
\begin{equation}  \label{cpilkku}
C' = 12^{-2} (C'')^{-2} ,
\qquad
C'' =  (1-a)^{-3/2} \Gamma(3/2) (\log_e2)^{-3/2} ,
\end{equation}
%where $0\leq a<1$ and $c>0$ are the constants giving the functional
%form of $\varrho(Q,{\cal G}_{\delta})$.
we have
$$
P \left( \sup_{g \in {\cal G}} \nu_n(Qg) \geq \xi \right)
\leq
\exp\left\{ -\frac{n\xi^2C'}{c^2 R^{2-2a}} \right\} .
$$
%where $\nu_n$ is the centered empirical process defined in (\ref{nun-def}).
\end{lemma}

{\em Proof.}
The proof uses the chaining technique.
The chaining technique was developed by Kolmogorov.
An analogous lemma in the direct case is for example
Lemma~3.2 in \cite{vanGe00}.
The basic difference to the direct case is visible in eq.~(\ref{basic-diffe}).
Let us denote $R_k=2^{-k}R$,
$
N_k = \#{\cal G}_{R_k}
$
and $H_k = \log_e N_k$, where $k=0,1,\ldots$.
For each $g \in {\cal G}$, let $h_g^k$ be a member of $R_k$
covering set of ${\cal G}$ such that
$\|g-h_g^k\|_{2} \leq R_k$.
We may write every $g \in {\cal G}$ with telescoping as
$$
g = \sum_{k=1}^{\infty} \left( h_g^k - h_g^{k-1} \right)
$$
where $h_g^0 \equiv 0$ and the convergence is in $L_2$.
Let $\eta_k>0$ be such that $\sum_{k=1}^{\infty} \eta_k \leq 1$.
We will define $\eta_k$ in (\ref{etakdef-white}).
Then
\begin{equation} %{eqnarray}
%\lefteqn{
\label{1stdiscompo-white}
P \left( \sup_{g \in {\cal G}} \nu_n(Qg) \geq \xi \right)
%} \\
%&&  \nonumber
\leq
\sum_{k=1}^{\infty}
P \left( \sup_{g \in {\cal G}} \nu_n\left(Q(h_g^k - h_g^{k-1})\right)
\geq \xi \eta_k \right) .
\end{equation} %{eqnarray}
We have
$$
\# \left\{  h_g^k - h_g^{k-1} :  g\in {\cal G}  \right\}
\leq
N_k N_{k-1}  \leq N_k^2  .
$$
We have
\begin{eqnarray}
\max\left\{ \left\| Q(h_g^k - h_g^{k-1}) \right\|_{2}
:  g \in {\cal G} \right\}
& \leq &   \nonumber
T_k \max\left\{ \left\| h_g^k - h_g^{k-1} \right\|_{2}
:  g \in {\cal G} \right\}
\\ & \leq &  \label{basic-diffe}
3 T_k R_k ,
\end{eqnarray}
where we denote $T_k=\varrho(Q,{\cal G}_{R_k})$,
when $\varrho(Q,{\cal G}_{\delta})$ is defined in (\ref{t-dense}),
and we used the fact
$$
\| h_g^k - h_g^{k-1} \|_{2}
\leq
\| h_g^k - g \|_{2}
+
\| h_g^{k-1}  -g \|_{2}
\leq
2^{-k}R + 2^{-k+1}R
=
3 R_k .
$$
When $W\sim N(0,\sigma^2)$, $\xi>0$, then
$
P(W>\xi)  \leq 2^{-1} \exp\{-\xi^2/(2\sigma^2)\} ,
$
see for example \cite{Dudle99}, Proposition 2.2.1.
We have that
$\nu_n(Q(h_g^k - h_g^{k-1})) \sim N(0,n^{-1}\| Q(h_g^k - h_g^{k-1}) \|_{2}^2)$
and thus
\begin{equation}  \label{hoefappli-white}
P \left( \sup_{g \in {\cal G}} \nu_n\left(Q(h_g^k - h_g^{k-1})\right)
\geq \xi \eta_k \right)
\leq
N_k^2 2^{-1} \exp\left\{ - \,\frac{1}{2}\,
\frac{n \xi^2 \eta_k^2}{3^2T_k^2R_k^2}  \right\} .
\end{equation}
Now we choose
\begin{equation}  \label{etakdef-white}
\eta_k
=
3T_kR_k \max\left\{
\frac{8^{1/2} H_k^{1/2}}{n^{1/2} \xi}\, ,
%\,\frac{2^{-k(1-a)} k^{1/2}}{2\sqrt{2}}
c^{-1} R^{a-1} (C'k)^{1/2} 2
\right\} ,
\end{equation}
where %$a$ is defined in Assumption~2 of Theorem~\ref{chain-white} and
$C'$ is defined in (\ref{cpilkku}).
Then we may apply (\ref{hoefappli-white})
to continue (\ref{1stdiscompo-white}) with an upper bound
\begin{eqnarray}
\frac{1}{2} \sum_{k=1}^{\infty}
\exp\left\{  2H_k - \,\frac{1}{2}\,
\frac{n \xi^2 \eta_k^2}{3^2T_k^2R_k^2}  \right\}
& \leq &  \label{askel1-white}
\frac{1}{2} \sum_{k=1}^{\infty}
\exp\left\{  - \, \frac{1}{4}\,
\frac{n \xi^2 \eta_k^2}{3^2T_k^2R_k^2}  \right\}
\\
& \leq &  \label{askel2-white}
\frac{1}{2} \sum_{k=1}^{\infty}
\exp\left\{  - \,\frac{n \xi^2 C' k}{c^2 R^{2-2a}} \right\}
\\
& \leq &   \label{final-white}
\exp\left\{ - \,\frac{n \xi^2 C'}{c^2 R^{2-2a}} \right\}  .
\end{eqnarray}
In (\ref{askel1-white}) we applied (\ref{etakdef-white}), which implies
$2H_k \leq n\xi^2\eta_k^2/(4\cdot 3^2 T_k^2R_k^2)$,
when we apply the first term in the maximum.
In (\ref{askel2-white}) we applied also (\ref{etakdef-white}), which implies
$
\eta_k^2/(4\cdot 3^2T_k^2R_k^2)
\geq
C'k /[c^2R^{2-2a}]
$
where we applied the second term in the maximum.
In (\ref{final-white}) we applied that for
$0 \leq b \leq 1/2$,
$\sum_{k=1}^{\infty} b^k = b/(1-b) \leq 2b$.
Here we need that
$\exp\left\{ - n \xi^2 C'/ [c^2R^{2-2a}]\right\} \leq 1/2$,
that is,
$
\xi \geq c R^{1-a} \left( \frac{\log_e 2}{nC'} \right)^{1/2}
$
which is implied by (\ref{xicond-white}).
We need to check that $\sum_{k=1}^{\infty} \eta_k \leq 1$.
Since 
$\delta \mapsto \varrho(Q,{\cal G}_{\delta}) \sqrt{\log_e(\#{\cal G}_{\delta})}$
is decreasing,
\begin{equation}  \label{kaneetti1}
\sum_{k=1}^{\infty} T_k R_k H_k^{1/2}
=
2\sum_{k=1}^{\infty} 2^{-k-1}R
T_{2^{-k}R} \sqrt{\log_e(\#{\cal G}_{2^{-k}R})}
%\sum_{k=1}^{\infty} (c R^{-a}2^{ak}) \cdot (3R2^{-k}) \cdot  H_k^{1/2}
%3c R^{1-a}\sum_{k=1}^{\infty} 2^{(a-1)k}  H_k^{1/2}
%
%\int_0^R T_u  \sqrt{\log_e(\#{\cal G}_{u})} \, du
\leq
2 G(R)  .
\end{equation}
We apply the assumption that
$T_k = T_{2^{-k}R} = cR^{-a} 2^{ak}$ to get
\begin{eqnarray}
\sum_{k=1}^{\infty} k^{1/2} T_kR_k
& = &     \nonumber
cR^{1-a} \sum_{k=1}^{\infty} k^{1/2} 2^{-(1-a)k}
\\ & = &
cR^{1-a} \lim_{K\to\infty} K^{3/2} \int_0^1 t^{1/2} 2^{-(1-a)Kt} \, dt
\\ & = &     \nonumber
cR^{1-a} (1-a)^{-3/2} \int_0^{\infty} u^{1/2} 2^{-u} \,du
\\ & = &  \label{kaneetti2}
c R^{1-a} C'',
\end{eqnarray}
where $C''$ is defined in (\ref{cpilkku}).
We have from (\ref{kaneetti1}) and  (\ref{kaneetti2}) that
$$
\sum_{k=1}^{\infty} \eta_k
\leq
\frac{8^{1/2} 6G(R)}{n^{1/2}\xi}
\, + \,
6 \sqrt{C'} C''
\leq
\frac{1}{2} + \frac{1}{2}
=
1 ,
$$
when $\xi \geq 28^{1/2}6 G(R) n^{-1/2}$,
which is guaranteed by (\ref{xicond-white}),
and $C'$ is chosen as in (\ref{cpilkku}).
The lemma follows from
(\ref{1stdiscompo-white}),
(\ref{hoefappli-white}),
and (\ref{final-white}).
\hspace*{\fill}   $\Box$  %$\square$ %$\Box$

%a<-0.9
%x<-seq(1,10,0.001)
%y<-x^(1/2)*2^(-(1-a)*x)
%plot(x,y)

\subsection{Density estimation}  \label{appe-empiri}

Lemma \ref{telehoef-density}
gives an exponential bound for the tail probability
in the case of density estimation.

\begin{lemma}  \label{telehoef-density}
Let $Y_1,\ldots ,Y_n\in {\bf R}^d$ be i.i.d.~with density $Af$,
and let the centered empirical process $\nu_n$ be defined in (\ref{nun-def}).
Assume that $\|Af\|_{\infty} \leq B_{\infty}$.
Let ${\cal G}\subset L_2({\bf R}^d)$
be such that
$\sup_{g\in{\cal G}} \|g\|_2 \leq R$.
Denote with ${\cal G}_{\delta}$ a $\delta$-bracketing net of ${\cal G}$,
$\delta>0$.
Denote
${\cal G}_{\delta}^L = \{ g^L : (g^L,g^U) \in {\cal G}_{\delta}\}$
and
${\cal G}_{\delta}^U = \{ g^U : (g^L,g^U) \in {\cal G}_{\delta}\}$.
Assume that
$\sup_{g\in {\cal G}_R^L \cup {\cal G}_R^U} \|Qg\|_{\infty} \leq B_{\infty}'$.
Assume that $\delta \mapsto \varrho_{den}(Q,{\cal G}_{\delta})
\sqrt{\log_e(\#{\cal G}_{\delta})}$
is decreasing on $(0,R]$,
where $\varrho_{den}(Q,{\cal G}_{\delta})$ is defined in (\ref{td-dens})
and
assume that the entropy integral $G(R)$
defined in (\ref{ent-int-density}) is finite.
Assume that $\varrho_{den}(Q,{\cal G}_{\delta}) = c\delta^{-a}$,
where $0\leq a <1$ and $c>0$.
Then for all
\begin{equation}  \label{xicond-density}
\xi \geq  n^{-1/2}\, \tilde{G}(R),
\end{equation}
where
\begin{eqnarray}   \nonumber
\tilde{G}(R)
& = &
B_{\infty}^{1/2} (9^2+96\cdot2^{-2a})^{1/2}
\\ && \times   \label{xicond-density-lisa}
\max\left\{ 24\sqrt{2}   G(R) ,
\,\, 4(\log_e(2))^{-1} (1-a)^{-3/2} \Gamma(3/2)   c R^{1-a}  \right\} ,
\end{eqnarray}
%\begin{equation}  \label{xicond-density-lisa}
%\tilde{G}(R)
%=
%B_{\infty}^{1/2} (9^2+96\cdot2^{-2a})^{1/2}
%\max\left\{ 24\sqrt{2}   G(R) ,
%\,\, 4 \sqrt{\log_e 2} C''  c R^{1-a}  \right\} ,
%\max\left\{ 24\sqrt{2} (9^2+96\cdot2^{-2a})^{1/2}  G(R) ,
%\,\, \sqrt{\log_e 2 / C'}  c R^{1-a}  \right\} ,
%\end{equation}
%$$ %\begin{equation}  \label{cpilkku-dens}
%C' = 4^{-2} ( C'')^{-2} (9^2+96\cdot2^{-2a})^{-1}  ,
%\qquad
%C'' = (1-a)^{-3/2} \Gamma(3/2) (\log_e2)^{-3/2} ,
%$$ %\end{equation}
%where $0\leq a<1$ and $c>0$ are the constants giving the functional
%form of $\varrho(Q,{\cal G}_{\delta})$.
we have
\begin{eqnarray*}
\lefteqn{
P \left( \sup_{g \in {\cal G}} \nu_n(Qg) \geq \xi \right)
} \\ & \leq &
4 \exp\left\{ -\frac{n\xi^2C'}{B_{\infty} c^2 R^{2-2a}} \right\}
+
2 \#{\cal G}_R \exp\left\{ -\, \frac{1}{12} \,
\frac{ n \xi^2}{B_{\infty}c^2R^{2(1-a)} + 2\xi B_{\infty}'/9} \right\} ,
\end{eqnarray*}
where $\nu_n$ is the centered empirical process defined in (\ref{nun-def}).
\end{lemma}

% a<-seq(0,1,0.01)
% cp<-(1-a)^(-3/2)*gamma(3/2)*log(2)^(-3/2)
% y<-4^(-2)*cp^(-2)*(9^2+96*2^(-2*a))^(-1)
% plot(a,y)

{\em Proof.}
We use the chaining technique with truncation.
The basic difference to the direct case is visible in
(\ref{basic-diffe-density}) and
(\ref{basic-diffe-density-0}).
The technique was used in the direct case by
\cite{Bass85},
\cite{Ossia87},
\cite{BirgeMassa93}, Proposition~3,
\cite{vanGe00}, Theorem~8.13.
Let us denote
$R_k=2^{-k}R$,
$
N_k =\#{\cal G}_{R_k}
$
and $H_k = \log_e N_k$, for $k=0,1,\ldots$.
Let us denote
%${\cal G}_k={\cal G}_{R_k}$
%a $R_k$-bracketing net of ${\cal G}$,
%and
$T_k=\varrho_{den}(Q,{\cal G}_{R_k})$, where 
$\varrho_{den}(Q,{\cal G}_{\delta})$
is defined in (\ref{td-dens}).
For each $g \in {\cal G}$, let $(h_g^{k,L},h_g^{k,U})$
be the member of the bracketing net ${\cal G}_{R_k}$,
such that $h_g^{k,L} \leq g \leq h_g^{k,U}$.
We may write every $g \in {\cal G}$ with telescoping as
$$
g =
g - h_g^{\kappa_g,L}
+ \sum_{k=1}^{\kappa_g} \left( h_g^{k,L} - h_g^{k-1,L} \right)
+ h_g^{0,L} ,
$$
where
$$
\kappa_g
=
\left\{ \begin{array}{ll}
\min\left\{ 0 \leq k \leq K-1 : Q\Delta_g^k \geq \beta_k \right\},
&
\mbox{ if } Q\Delta_g^k \geq \beta_k \mbox{ for some } 0\leq k \leq K-1
\\
K,
& \mbox{ otherwise, }
\end{array}  \right.
$$
where $K\geq 1$ is defined in (\ref{bigk}),
$$
\Delta_g^k = h_g^{k,U} - h_g^{k,L} ,
%\Delta_g^k = \tilde{h}_g^{k,U} - \tilde{h}_g^{k,L} ,
%\qquad
%\tilde{h}_g^{k,U}(x) = \min_{0 \leq l \leq k} h_g^{l,U}(x),
%\qquad
%\tilde{h}_g^{k,L}(x) = \max_{0 \leq l \leq k} h_g^{l,L}(x) ,
$$
%$x\in {\bf R}^d$
and
\begin{equation}  \label{betak-bound}
\beta_k = \frac{12B_{\infty} T_k^2R_k^2}{\xi}  \, .
\end{equation}
Then,
\begin{eqnarray}
P \left( \sup_{g \in {\cal G}} \nu_n(Qg) \geq \xi \right)
& \leq &     \nonumber
P \left( \sup_{g \in {\cal G}} \sum_{k=1}^{\kappa_g}
\nu_n\left(Q( h_g^{k,L} - h_g^{k-1,L})\right)
\geq \xi/3  \right)
\\ && +       \nonumber
P \left( \sup_{g \in {\cal G}} \nu_n\left( Q( g - h_g^{\kappa_g,L}) \right)
\geq \xi/3 \right)
\\ && +     \nonumber
P \left( \sup_{g \in {\cal G}} \nu_n\left( Qh_g^{0,L} \right)
\geq \xi/3 \right)
\\ & \stackrel{def}{=} & \label{1stdiscompo-density}
P_I + P_{II} + P_{III} .
\end{eqnarray}

\bigskip \noindent {\em Term $P_{I}$.}
We have
\begin{eqnarray*}
\sup_{g \in{\cal G}} \sum_{k=1}^{\kappa_g}
\nu_n\left( Q(h_g^{k,L} - h_g^{k-1,L}) \right)
& = &
\sup_{g \in{\cal G}} \sum_{k=1}^{K}
I_{\{1,\ldots ,\kappa_g\}}(k) \nu_n\left( Q(h_g^{k,L} - h_g^{k-1,L})\right)
\\ & \leq &
\sum_{k=1}^{K}  \sup_{g \in{\cal G}}
I_{\{1,\ldots ,\kappa_g\}}(k) \nu_n\left(Q( h_g^{k,L} - h_g^{k-1,L})\right) .
\end{eqnarray*}
Let us denote
\begin{equation}  \label{etakdef-density}
\eta_k
=
(9^2+96\cdot 2^{-2a})^{1/2}  T_kR_k
\max\left\{
\frac{8^{1/2} B_{\infty}^{1/2} H_k^{1/2}  }{n^{1/2} \xi}\,
,  c^{-1} R^{a-1} (C'k)^{1/2} 2
\right\} ,
\end{equation}
where $C'$ is defined by %in (\ref{cpilkku-dens}).
\begin{equation}  \label{cpilkku-dens}
C' = 4^{-2} ( C'')^{-2} (9^2+96\cdot2^{-2a})^{-1}  ,
\qquad
C'' = (1-a)^{-3/2} \Gamma(3/2) (\log_e2)^{-3/2} .
\end{equation}
We have defined $\eta_k$ in (\ref{etakdef-density}) so that
$\eta_k>0$ and $\sum_{k=1}^{\infty} \eta_k \leq 1$, which is proved in
(\ref{etaksum}).
Then, %for $\eta_k$ defined in (\ref{etakdef-density}),
\begin{equation} \label{pii-aku}
P_{I}
\leq
\sum_{k=1}^{K}
P \left( \sup_{g \in {\cal G}}
I_{\{1,\ldots ,\kappa_g\}}(k) \nu_n\left(Q( h_g^{k,L} - h_g^{k-1,L})\right)
\geq \eta_k\xi/3  \right)  .
\end{equation}
We have
\begin{equation}  \label{cardibound-3}
\# \left\{  h_g^{k,L} - h_g^{k-1,L} :  g\in {\cal G}  \right\}
\leq
N_k N_{k-1}  \leq N_k^2  .
\end{equation}
Also,
\begin{eqnarray}
\lefteqn{ \nonumber
\max\left\{ E\left| Q(h_g^{k,L} - h_g^{k-1,L}) \right|^{2}
:  g \in {\cal G} \right\}
} \\ & \leq &   \nonumber
B_{\infty} \max\left\{ \left\| Q(h_g^{k,L} - h_g^{k-1,L}) \right\|_{2}^2
:  g \in {\cal G} \right\}
\\ & \leq &   \nonumber
B_{\infty} T_k^2 \max\left\{ \left\| h_g^{k,L} - h_g^{k-1,L} \right\|_{2}^2
:  g \in {\cal G} \right\}
\\ & \leq &  \label{basic-diffe-density}
B_{\infty} 3^2 T_k^2 R_k^2 ,
\end{eqnarray}
because
$$
\left\| h_g^{k,L} - h_g^{k-1,L} \right\|_{2}
\leq
\left\| h_g^{k,L} - g \right\|_{2}
+
\left\| h_g^{k-1,L}  -g \right\|_{2}
\leq
2^{-k}R + 2^{-k+1}R
=
3 R_k .
$$
When $k \leq \kappa_g$, then
$$
Q(h_g^{k,L} - h_g^{k-1,L}) \leq Q\Delta_g^{k-1} \leq \beta_{k-1} ,
$$
which implies
\begin{equation}  \label{unibound-3}
\left| Q(h_g^{k,L} - h_g^{k-1,L}) - E Q(h_g^{k,L} - h_g^{k-1,L}) \right|
\leq 2\beta_{k-1} .
\end{equation}
Thus, applying
(\ref{cardibound-3}),
(\ref{basic-diffe-density}),
(\ref{unibound-3}),
by Bernstein's inequality,
\begin{eqnarray}
\lefteqn{  \nonumber
P \left( \sup_{g \in {\cal G}} I_{\{1,\ldots ,\kappa_g\}}(k)
\nu_n\left( Q(h_g^{k,L} - h_g^{k-1,L})\right)
\geq \xi \eta_k/3 \right)
} \\ & \leq &  \nonumber %\label{hoefappli-density}
N_k^2 \exp\left\{ - \,\frac{1}{2}\,
\frac{n (\xi\eta_k/3)^2}{3^2B_{\infty}T_k^2R_k^2+
2\beta_{k-1}\xi\eta_k/9}  \right\}
\\ & \leq &        \label{drop}
\exp\left\{  2H_k
 - \,\frac{1}{2}\,
\frac{n (\xi\eta_k)^2}{3^2(3^2+24\cdot2^{2(1-a)}/9)B_{\infty}T_k^2R_k^2}
\right\}
\\ & \leq &        \label{drop2}
\exp\left\{
 - \,\frac{1}{4}\,
\frac{n (\xi\eta_k)^2}{(9^2+96\cdot2^{-2a})B_{\infty}T_k^2R_k^2}
\right\}
\\ & \leq &  \label{askel2-density}
\exp\left\{  - \,\frac{n \xi^2 C' k}{c^2B_{\infty} R^{2-2a}}  \right\} .
\end{eqnarray}
In (\ref{drop}) we applied the fact
$\beta_{k-1}\xi\eta_k \leq 12B_{\infty}2^{2(1-a)} T_k^2R_k^2$,
which follows since
$T_kR_k = cR_k^{1-a} = 2^{1-a}T_{k+1}R_{k+1}$, which implies
\begin{equation}  \label{beta-yla}
\beta_k
\leq
\frac{12B_{\infty} 2^{2(1-a)} T_{k+1}^2R_{k+1}^2}{\eta_{k+1} \xi} \, ,
\end{equation}
since $0<\eta_k\leq 1$,
where $\eta_k$ is defined in (\ref{etakdef-density}).
In (\ref{drop2}) we applied the first term in the maximum in
(\ref{etakdef-density}) which implies
$2H_k \leq 4^{-1}n(\xi\eta_k)^2/[(9^2+96\cdot2^{-2a})
B_{\infty} T_k^2R_k^2]$.
In (\ref{askel2-density}) we applied the second term in the maximum in
(\ref{etakdef-density}),
which implies
$
\eta_k^2/[4\cdot (9^2+96\cdot2^{-2a}) T_k^2R_k^2]
\geq
C'k / [c^2R^{2-2a}] .
$
We may continue (\ref{pii-aku}) with an upper bound
\begin{equation}
\sum_{k=1}^{\infty}
\exp\left\{  - \,\frac{n \xi^2 C' k}{c^2B_{\infty} R^{2-2a}}  \right\}
\leq    \label{final-density}
2 \exp\left\{ - \,\frac{n \xi^2 C'}{c^2B_{\infty} R^{2-2a}}  \right\}   .
\end{equation}
%In (\ref{final-density})
We applied the fact that for
$0 \leq a \leq 1/2$,
$\sum_{k=1}^{\infty} a^k = a/(1-a) \leq 2a$.
Here we need that
$\exp\left\{ - n \xi^2 C'/[c^2B_{\infty}R^{2-2a}] \right\} \leq 1/2$,
that is we need,
$
\xi \geq \left( \frac{\log_e 2}{nC'} \right)^{1/2}
c B_{\infty}^{1/2} R^{1-a}
$
which is implied by (\ref{xicond-density}).

% 3^2*(3^2+24*2^2/9) =177

% 6^2*(1+2^2/3) =84

\bigskip \noindent {\em Term $P_{II}$.}
We have
$$
g - h_g^{k,L} \leq \Delta_g^{k}
$$
and thus
$$
\nu_n\left(Q(  g - h_g^{\kappa_g,L}) \right)
\leq
\nu_n\left( Q\Delta_g^{\kappa_g} \right )
+ 2 E\left| Q\Delta_g^{\kappa_g} \right| .
$$
Here we used the assumption that operator $Q$ preserves positivity
($g \geq 0$ implies $Qg \geq 0$).
We have for $k=0,\ldots ,K$,
\begin{eqnarray}
\max\left\{ E \left| Q\Delta_g^{k} \right|^2 : g \in {\cal G} \right\}
& \leq &   \nonumber
B_{\infty}
\max\left\{ \left\| Q\Delta_g^{k} \right\|_{2}^2 : g \in {\cal G} \right\}
\\ & \leq &   \nonumber
B_{\infty} T_k^2 \max\left\{ \left\| \Delta_g^{k} \right\|_{2}^2 :
g \in {\cal G} \right\}
\\ & \leq &  \label{basic-diffe-density-0}
B_{\infty} T_{k}^2 R_k^2 .
\end{eqnarray}
When $\kappa_g=k$, then $Q\Delta_g^{\kappa_g} \geq \beta_k$,
for $k=0,\ldots, K-1$.
Thus, for $\kappa_g=k$, $k=0,\ldots, K-1$,
using (\ref{betak-bound}),
$$
E| Q\Delta_g^{\kappa_g}|
\leq
\beta_k^{-1} E| Q\Delta_g^{k}|^2
%\leq
%\beta_k^{-1} B_{\infty} \| Q\Delta_g^{k}\|_2^2
\leq
\beta_k^{-1} B_{\infty} T_k^2 R_k^2
\leq
\xi/12 ,
$$
and for $\kappa_g=K$,
$$
E| Q\Delta_g^{\kappa_g}|
\leq
\left( E| Q\Delta_g^{K}|^2 \right)^{1/2}
\leq
B_{\infty}^{1/2} T_KR_K
\leq
\xi/12 ,
$$
when we choose
\begin{equation}  \label{bigk}
K =  \min\left\{ k \geq 1 : 12 B_{\infty}^{1/2} T_kR_k  < \xi  \right\}.
\end{equation}
Thus
$$
P\left( \sup_{g\in {\cal G}} 2 E| Q\Delta_g^{\kappa_g}| > \xi/6\right) = 0 .
$$
Define
$$
{\cal G}^{(k)} = \{ g \in {\cal G}  : \kappa_g = k \},
\qquad
k=0,\ldots , K ,
$$
so that
$
{\cal G} = \bigcup_{k=0}^K {\cal G}^{(k)} .
$
%We have that
%$$
%\sup_{g\in {\cal G}} \nu_n\left( Q(g-h_g^{\kappa_g,L}) \right)
%\leq
%\sum_{k=1}^K \sup_{g\in {\cal G}^{(k)}}
%\nu_n\left( Q(g-h_g^{\kappa_g,L}) \right)
%\leq
%\sum_{k=1}^K \sup_{g\in {\cal G}^{(k)}}
%\nu_n\left( Q\Delta_g^{k} \right)   .
%$$
Then,
\begin{eqnarray} \nonumber
P_{II}
& \leq &
P\left( \sup_{g\in {\cal G}}
\nu_n\left( Q\Delta_g^{\kappa_g} \right) \geq \xi/6 \right)
\leq
\sum_{k=0}^{K}
P \left( \sup_{g \in {\cal G}^{(k)}} \nu_n\left(Q\Delta_g^{k}\right)
\geq \xi/6 \right) \nonumber
\\ & = &
P_{II}^{(0)} + P_{II}^{(1)} ,
\label{noolla-density}
\end{eqnarray}
where
$$
P_{II}^{(0)}
=
P \left( \sup_{g \in {\cal G}^{(0)}} \nu_n\left(Q\Delta_g^{0}\right)
\geq \xi/6 \right) ,
\,\,
P_{II}^{(1)}
=
\sum_{k=1}^{K}
P \left( \sup_{g \in {\cal G}^{(k)}} \nu_n\left(Q\Delta_g^{k}\right)
\geq \xi/6 \right) .
$$
We have
\begin{equation}  \label{cardi-for}
\# \left\{  \Delta_g^{k}  :  g\in {\cal G}^{(k)}  \right\}
\leq
\# \left\{  \Delta_g^{k}  :  g\in {\cal G}  \right\}
\leq
N_{k} .
\end{equation}
It holds that
\begin{equation}  \label{unifor-for-pikku}
\left| Q\Delta_g^{0} - E Q\Delta_g^{0} \right| \leq 4B_{\infty}' .
\end{equation}
We have, using
(\ref{basic-diffe-density-0}),
(\ref{cardi-for}),
(\ref{unifor-for-pikku}),
by Bernstein's inequality,
\begin{equation}     \label{drop2-fin-not}
P_{II}^{(0)}
\leq
N_0 \exp\left\{ - \,\frac{1}{2}\,
\frac{n (\xi/6)^2}{B_{\infty} T_0^2R_0^2
+ 2B_{\infty}'\xi/9}  \right\} .
\end{equation}
Let us turn to $P_{II}^{(1)}$.
For $\kappa_g =k$ (that is, when $g\in {\cal G}^{(k)}$),
for $k=1,\ldots ,K$,
$$
Q\Delta_g^{k} \leq Q\Delta_g^{k-1} \leq \beta_{k-1} ,
$$
which implies
\begin{equation}  \label{unifor-for}
\left| Q\Delta_g^{k} - E Q\Delta_g^{k} \right| \leq 2 \beta_{k-1} .
\end{equation}
Thus, using
(\ref{basic-diffe-density-0}),
(\ref{cardi-for}),
(\ref{unifor-for}),
the fact that $0<\eta_k\leq 1$,
where $\eta_k$ is defined in (\ref{etakdef-density}),
by Bernstein's inequality,
for $k=1,\ldots ,K$,
\begin{eqnarray}
P \left( \sup_{g \in {\cal G}^{(k)}} \nu_n( Q\Delta_g^{k})
\geq \xi/6 \right)
& \leq &  \nonumber
P \left( \sup_{g \in {\cal G}^{(k)}} \nu_n( Q\Delta_g^{k})
\geq \xi\eta_k/6 \right)
\\ & \leq &  \nonumber
N_k \exp\left\{ - \,\frac{1}{2}\,
\frac{n (\xi\eta_k/6)^2}{B_{\infty}T_k^2R_k^2
+ \beta_{k-1}\xi\eta_k/9}  \right\}
\\ & \leq &        \label{drop-0}
\exp\left\{  H_k
 - \,\frac{1}{2}\,
\frac{n (\xi\eta_k)^2}{6^2(1+2^{2(2-a)}/3) B_{\infty} T_k^2R_k^2}
\right\}
\\ & \leq &        \label{drop2-0}
\exp\left\{
 - \,\frac{1}{4}\,
\frac{n (\xi\eta_k)^2}{(6^2+48\cdot2^{-2a})B_{\infty} T_k^2R_k^2}
\right\}
\\ & \leq &   \label{askel2-density-0}
\sum_{k=1}^{\infty}
\exp\left\{  - \,\frac{n \xi^2 C' k}{B_{\infty}c^2R^{2-2a}}  \right\} .
\end{eqnarray}
In (\ref{drop-0}) we applied the fact
$\beta_{k-1}\xi\eta_k \leq 12B_{\infty}2^{2(1-a)}T_k^2R_k^2$,
which follows by using (\ref{beta-yla}).
In (\ref{drop2-0}) we applied the first term in the maximum in
(\ref{etakdef-density}) which implies
$H_k \leq 4^{-1}n(\xi\eta_k)^2/[(6^2+48\cdot2^{-2a}) B_{\infty}T_k^2R_k^2]$,
since $2^{-1}(6^2+48\cdot 2^{-2a}) \leq 9^2+96\cdot 2^{-2a}$.
%We may continue (\ref{noolla-density}) with an upper bound
In (\ref{askel2-density-0}) we applied the second term in the maximum in
(\ref{etakdef-density}),
which implies
$
\eta_k^2/[4\cdot (6^2+48\cdot2^{-2a}) T_k^2R_k^2]
\geq
C'k / [c^2R^{2-2a}] .
$
We get
\begin{equation}
P_{II}^{(1)}
\leq
\sum_{k=1}^{\infty}
\exp\left\{  - \,\frac{n \xi^2 C' k}{B_{\infty}c^2R^{2-2a}}  \right\}
\leq    \label{pi-end}  %\label{final-density-0}
2 \exp\left\{ - \,\frac{n \xi^2 C'}{B_{\infty}c^2R^{2-2a}}  \right\}   .
\end{equation}
In (\ref{pi-end}) we applied that for
$0 \leq a \leq 1/2$,
$\sum_{k=1}^{\infty} a^k = a/(1-a) \leq 2a$.
Here we need that
$\exp\left\{ - n \xi^2 C'/ [B_{\infty} c^2R^{2-2a}] \right\} \leq 1/2$,
that is we need,
$
\xi \geq \left( \frac{\log_e 2}{nC'} \right)^{1/2} B_{\infty} cR^{1-a}
$
which is implied by (\ref{xicond-density}).

\bigskip \noindent {\em Term $P_{III}$.}
We have first,
$$
\# \left\{ h_g^{0,L} : g \in {\cal G} \right\}
\leq
N_0 ,
$$
second
$$
\sup_{g \in {\cal G}} E \left| Q h_g^{0,L} \right|^2
\leq
B_{\infty} \sup_{g \in {\cal G}} \left\| Q h_g^{0,L} \right\|_2^2
\leq
B_{\infty} T_0^2R_0^2
=
B_{\infty} c^2R^{2-2a} ,
$$
and third
$$
\sup_{g \in {\cal G}} \left\| Q h_g^{0,L} \right\|_{\infty}
\leq
B_{\infty}' .
$$
Thus, by Bernstein's inequality
\begin{equation}   \label{drop2-fin}
P_{III}
\leq
N_0 \exp\left\{ -\, \frac{1}{2} \,
\frac{ n (\xi/3)^2}{B_{\infty}T_0^2R_0^2 + \xi 2B_{\infty}'/9} \right\} .
\end{equation}

\bigskip \noindent {\em Finishing the proof.}
The lemma follows from
(\ref{1stdiscompo-density}),
(\ref{final-density}),
(\ref{drop2-fin-not}),
%(\ref{noolla-density}),
(\ref{pi-end}),
and
(\ref{drop2-fin}),
after checking some final facts.
We need to check that $\sum_{k=1}^{\infty} \eta_k \leq 1$.
Applying the calculations in (\ref{kaneetti1}) and (\ref{kaneetti2})
we get
\begin{equation}  \label{etaksum}
\sum_{k=1}^{\infty} \eta_k
\leq
(9^2+96\cdot 2^{-2a})^{1/2}
\left(
\frac{8^{1/2}B_{\infty}^{1/2} 2G(R)}{n^{1/2}\xi}
\, + \,
2\sqrt{C'} C''  \right)
\leq
\frac{1}{2} + \frac{1}{2}
=
1 ,
\end{equation}
when
$\xi \geq 2\cdot 8^{1/2}6 G(R) n^{-1/2}
(9^2+96\cdot 2^{-2a})^{1/2} B_{\infty}^{1/2}$,
which is guaranteed by (\ref{xicond-density}),
and $C'$ is chosen as in (\ref{cpilkku-dens}).
\hspace*{\fill}   $\Box$  %$\square$ %$\Box$

\begin{remark}{\em
When in addition $\xi$ satisfies
\begin{equation}  \label{xi-yla}
2  \sqrt{44 \log_e(\#{\cal G}_R)} B_{\infty}^{1/2}c R^{1-a} n^{-1/2}
\leq
\xi
\leq  B_{\infty} c^2R^{2(1-a)}/B_{\infty}' ,
\end{equation}
then
$$ %\begin{equation}  \label{sievi}
\#{\cal G}_R \exp\left\{ -\, \frac{1}{12} \,
\frac{ n \xi^2}{B_{\infty}c^2R^{2(1-a)} + 2\xi B_{\infty}'/9} \right\}
\leq
\exp\left\{ -\,
\frac{ n \xi^2 C'}{B_{\infty}c^2R^{2(1-a)}} \right\} .
$$  %\end{equation}
Indeed, we may continue (\ref{drop2-fin-not}) by
\begin{eqnarray}
P_{II}^{(1)}  \nonumber
& \leq &
N_0 \exp\left\{ - \,\frac{1}{2}\,
\frac{n (\xi/6)^2}{B_{\infty} T_0^2R_0^2
+ 2B_{\infty}'\xi/9}  \right\}
\\ & \leq &      \label{piii-end-not}
\exp\left\{ H_0 -\, \frac{1}{2} \,
\frac{ n \xi^2}{6^2(1+2/9)B_{\infty}c^2R^{2(1-a)} } \right\}
\\ & \leq &      \label{drop2-fin-not-not}
\exp\left\{ -\, \frac{1}{4} \,
\frac{ n \xi^2}{44B_{\infty}c^2R^{2(1-a)}} \right\} .
\end{eqnarray}
In (\ref{piii-end-not}) we applied the upper bound in (\ref{xi-yla})
and the fact $T_0R_0=cR^{1-a}$.
In (\ref{drop2-fin-not-not}) we applied the lower bound in (\ref{xi-yla})
which implies the fact
$H_0 \leq 4^{-1}n\xi^2/[44 B_{\infty}c^2R^{2(1-a)}]$.
Also, we may continue (\ref{drop2-fin}) by
\begin{eqnarray}
P_{III}  \nonumber
& \leq &
N_0 \exp\left\{ -\, \frac{1}{2} \,
\frac{ n (\xi/3)^2}{B_{\infty}T_0^2R_0^2 + \xi 2B_{\infty}'/9} \right\}
\\ & \leq &      \label{piii-end}
\exp\left\{ H_0 -\, \frac{1}{2} \,
\frac{ n \xi^2}{3^2(1+2/9)B_{\infty}c^2R^{2-2a}} \right\}
\\ & \leq &      \label{drop2-fin-fin}
\exp\left\{ -\, \frac{1}{4} \,
\frac{ n \xi^2}{11 B_{\infty}c^2R^{2-2a}} \right\} .
\end{eqnarray}
In (\ref{piii-end}) where we applied the upper bound in (\ref{xi-yla}).
In (\ref{drop2-fin-fin}) we applied the lower bound in (\ref{xi-yla})
which implies
$H_0 \leq 4^{-1}n\xi^2/(11B_{\infty} c^2R^{2-2a})$.
}\end{remark}

\section{Introductory remarks}

We add a short introduction to the setting of the article,
in order to make the article more accessible to PhD students.

A quite general inverse problem could be described as a problem
where we want to recover a function $f:{\bf R}^d \to {\bf R}$
when we have only available some transform $Af$ of the
function. An important example is the sampling operator
$Af = (f(x_1),\ldots ,f(x_n)) \in {\bf R}^n$,
where
$x_1,\ldots ,x_n \in {\bf R}^d$ are fixed points.
Classical methods for recovering $f$ in this case 
include piecewise constant interpolation and various ways
to linearly interpolate the observed function values.
In statistics some kind of sampling operator is always involved
and thus recovering $f$ from noisy data
$f(x_i) + \epsilon_i$, $i=1,\ldots, n$,
where $\epsilon_i$ are error terms,
would not be called an inverse problem in statistics.
We mention three classical statistical inverse problems,
where function $f:{\bf R}^d \to {\bf R}$ has to be estimated
and $A$ is a fixed operator mapping functions 
${\bf R}^d \to {\bf R}$
to functions ${\bf Y} \to {\bf R}$,
where ${\bf Y}$ is some general space.

\begin{enumerate}
\item
{\em (Regression function estimation.)}
We observe data 
$$
Y_i 
=
(Af)(X_i) + \epsilon_i  \in {\bf R},
\qquad
i=1,\ldots ,n,
$$
where
$\epsilon_i \in {\bf R}$ are random errors
and 
$X_i \in {\bf Y}$ are random design points.

\item
{\em (Density estimation.)}
We observe identically distributed observations
$$
Y_1,\ldots , Y_n \in {\bf Y},
$$
whose common density is $Af$.

\item
{\em (Gaussian white noise model.)}
We observe a realization of the process
$$
dY_n(y) 
=
(Af)(y)\, dt + n^{-1/2} dW(y) ,
\qquad
y \in {\bf Y},
$$
where
$W(y)$ is a Wiener process on ${\bf Y}$.
When ${\bf Y} = {\bf R}$, then
we can define the process by 
$$
Y_n(y)
=
\int_{-\infty}^{y} (Af)(t) \, dt + n^{-1/2} W(y) ,
$$
where 
$W$ is the Brownian motion, or Wiener process, on the real line.
The Gaussian white noise model is rather close to the regression
function estimation when the error terms $\epsilon_i$ are Gaussian
and the design points $X_i$ are uniformly distributed in the unit square.
However, in the Gaussian white noise model we have eliminated
the problems related to interapolation since the function $Af$
is observed continuously and not in a finite number of design points.
Since the assumption of continuous observation is quite far from reality,
we can use inference in the Gaussian white noise model
only as a first approximation.
In addition, the assumption of the exact Gaussian distribution is 
very restrictive.
Due to the central limit theorem the Gaussian white noise
model is a relevant approximation also for the model of density estimation
and for the model of regression function estimation under
non-Gaussian noise.
\end{enumerate}

Let us now consider the estimation of a regression function 
(item 1 of the above list).
A common approach for regression function estimation is to
find the estimator $\tilde{f}$ as a solution of the minimization problem
\begin{equation}  \label{regressio}
\tilde{f}
=
\mbox{argmin}_{g \in {\cal F}}
\sum_{i=1}^n (Y_i - (Ag)(X_i))^2 ,
\end{equation}
where
${\cal F}$ is some class of functions ${\bf R}^d \to {\bf R}$.
Note that estimator $\hat{f}$ is a special case of 
the linear regression estimator 
$$
\tilde{f}(x) 
= 
\hat{\beta}_0 + \hat{\beta}_1^T x,
\qquad
(\hat{\beta}_0,\hat{\beta}_1)
=
\mbox{argmin}_{\beta_0 \in {\bf R}, \beta_1 \in {\bf R}^d}
\sum_{i=1}^n (Y_i - \beta_0 - \beta_1^T X_i)^2 ,
$$
when 
$A$ is the identity operator and
${\cal F} = 
\{ \beta_0 + \beta_1^Tx : \beta_0\in {\bf R},\beta_1 \in {\bf R}^d \}$. 
Estimator $\tilde{f}$, defined in (\ref{regressio}),
can be defined also as
\begin{equation}  \label{regressio-alt}
\tilde{f}
=
\mbox{argmin}_{g \in {\cal F}}
\left(
-\,\frac{2}{n} \sum_{i=1}^n Y_i\cdot  (Ag)(X_i) 
+
\frac{1}{n} \sum_{i=1}^n (Ag)^2(X_i) \right)  .
\end{equation}
The estimator which we have considered can be defined,
assuming now for simplicity that the design points $X_i$
have a known distribution $\nu$ on ${\bf Y}$,
\begin{equation}  \label{regressio-new}
\hat{f}
=
\mbox{argmin}_{g \in {\cal F}}
\left(
-\,\frac{2}{n} \sum_{i=1}^n Y_i\cdot  (Qg)(X_i) 
+
\int_{[0,1]^d} g^2 
 \right)  ,
\end{equation}
where
$Q = (A^{-1})^*$ is the adjoint of the inverse,
for the space $L_2(\nu)$.
Note that when $g : {\bf R}^d \to {\bf R}$,
then $Qg: {\bf Y} \to {\bf R}$.

When operator $B : H_1 \to H_2$ is defined as a mapping from a
Hilbert space $H_1$ to an another Hilbert space $H_2$,
then
the adjoint $B^*$ is defined as the operator
satisfying the equality
$$
\langle Bx ,y \rangle_2
=
\langle x , B^*y \rangle_1 ,
$$
where $\langle \cdot,\cdot \rangle_i$ are the inner products of the
Hibert spaces.
%(this generalizes the equality (\ref{q-def})
In the case when the Hilbert spaces are the Euclidean space: 
$H_1=H_2={\bf R}^d$, then the operators
are $d \times d$ real matrices, and we have
$\langle Bx ,y \rangle
=
\langle x , B^Ty \rangle$,
where
$B^T$ is the transpose of matrix $B$, and thus the adjoint
is equal to the transpose.
We have given further examples of adjoints in
(\ref{convo-adjoint}),
where the adjoint of the inverse of a convolution operator is given,
and in (\ref{defperusidea}),
where the adjoint of the inverse of the Radon transform is given.

The estimator defined by 
(\ref{regressio}) and (\ref{regressio-alt}) seems quite 
natural but we can justify the estimator defined
in (\ref{regressio-new}) by the following calculation.
We have, similarly as in (\ref{empi-moti}),
\begin{eqnarray*}
\| \hat{f}-f\|_2^2 - \| f\|_2^2
&=&  \nonumber
-2 \int_{{\bf R}^d} f\hat{f} + \|\hat{f}\|_2^2
\\ &=&  \nonumber
-2 \int_{{\bf Y}} (Af)(Q\hat{f}) \, d\nu + \|\hat{f}\|_2^2
\\ &\approx&  \nonumber
-\, \frac{2}{n} \sum_{i=1}^n Y_i \cdot (Q\hat{f})(X_i) 
+ \|\hat{f}\|_2^2 .
\end{eqnarray*}
The last approximation in the above calculation uses the fact
that the distribution of the design points is $\nu$.

\end{document}